\documentclass[10pt,reqno,a4paper]{article}

\usepackage[hypertex]{hyperref}
\usepackage{amsfonts}
\usepackage{latexsym}
\usepackage{amssymb,amsmath}
\usepackage{epsfig}

\newcommand{\C}{\mathbb C}
\newcommand{\R}{\mathbb R}
\newcommand{\N}{\mathbb N}
\newcommand{\Z}{\mathbb Z}

\newcommand{\T}{\mathbb T}
\newcommand{\TN}{\T ^n}

\newcommand{\be}{\begin{equation}}
\newcommand{\ee}{\end{equation}}
\newcommand{\ba}{\begin{eqnarray}}
\newcommand{\ea}{\end{eqnarray}}
\newcommand{\la}{\langle}
\newcommand{\ra}{\rangle}

\def\qd{{\rule{2.3mm}{2.3mm}}}
\def\qed{{\hfill$\quad$\qd\\}}

\def\ep{{\varepsilon}}

\newtheorem{thm}{Theorem}[section]
\newtheorem{rem}[thm]{Remark}
\newtheorem{lem}[thm]{Lemma}
\newtheorem{cor}[thm]{Corollary}
\newtheorem{prop}[thm]{Proposition}
\newtheorem{defi}[thm]{Definition}

\def\qd{\rule{2mm}{2mm}}

\begin{document}

\title{Control and Stabilization of the Nonlinear
Schr\"odinger Equation on Rectangles}

\author{
Lionel Rosier 
\thanks{Institut \'Elie Cartan,
UMR 7502 UHP/CNRS/INRIA,
B.P. 239,  54506 Vand\oe uvre-l\`es-Nancy Cedex, France.
({\tt rosier@iecn.u-nancy.fr}) }
\and 
Bing-Yu Zhang 
\thanks{Department of Mathematical Sciences,
University of Cincinnati, Cincinnati,
Ohio 45221, USA. 
({\tt bzhang@math.uc.edu})}
}

\maketitle



\begin{abstract}
This paper studies the local exact controllability and the local stabilization
of the semilinear Schr\"odinger equation posed on a product of $n$
intervals ($n\ge 1$).
Both internal and boundary controls are considered, and the results 
are given with periodic (resp. Dirichlet or Neumann) boundary
conditions. In the case of internal control, we obtain  
local controllability results which are sharp as far as the localization
of the control region and the smoothness of the state space are concerned.
It is also proved that for the linear Schr\"odinger equation with
Dirichlet control, the exact controllability holds in $H^{-1}(\Omega )$
whenever the control region contains a neighborhood of a vertex.



\end{abstract}
\section{Introduction}
\setcounter{equation}{0}

The control of the Schr\"odinger equation has received a lot of
attention in the last decades. (See e.g.
\cite{Zua-2} for an excellent review of the contributions up to 2003).  
Significant progresses have been made
for the linear Schr\"odinger equation on  its controllability and
stabilizability properties (see \cite{jaffard,KL,lebeau,liu,mach-2,
mach-3,miller,RTTT} for control issues, and \cite{BP,CG,CCG,mor,YY} for Carleman estimates and their applications to inverse problems).
For the control of the so-called {\em bilinear}
Schr\"odinger equation, in which the bilinear term is  linear in both 
the control and the state function, see e.g. 
\cite{BMS,CLP,BCG,Beauchard,BKP,MRT,ILT-2,BC,BS} and the references therein.

By contrast, the study of the nonlinear Schr\"odinger equation is still at
 its early stage. Recently,  Illner, Lange and  Teismann  
\cite{ILT-1,ILT-2} considered the internal controllability
 of the nonlinear Schr\"odinger equation posed on  a finite interval
 with periodic boundary conditions:
\be
\label{dim1}
iu_t+u_{xx}+f(u)=ia(x)h(x,t).
\ee
In (\ref{dim1}), $a$ denotes a smooth real function which is strictly
supported in $\T$, the one-dimensional torus.  They showed that the 
system \eqref{dim1}
 is locally exactly controllable in the space $H^1(\T )$.
 Their approach was based on the well-known Hilbert Uniqueness 
Method (HUM) and  Schauder's fixed point theorem.
 Later, Lange and  Teismann  \cite{LT-1} considered internal control for the
 nonlinear Schr\"odinger equation (\ref{dim1})  posed on a finite interval with
 the homogeneous Dirichlet boundary conditions 
\be
\label{dirichlet}
u(0,t)=u(\pi ,t)=0
\ee
and
 established  local exact controllability of the system
 (\ref{dim1})-(\ref{dirichlet}) in the space $H^1_0 (0, \pi)$  
around a special ground state of the
 system. Their approach was mainly based upon HUM and the implicit
 function theorem. Dehman,  G\'erard and  Lebeau 
\cite{DGL} studied  the internal control and stabilization of a class
 of defocusing  nonlinear Schr\"odinger equations
 posed on a two-dimensional compact Riemannian manifold  $M$ 
without boundary
$$
iu_t+\Delta u + f(u)=ia(x)h(x,t).
$$
They demonstrated, in particular, that the system is (semiglobally)
exactly controllable and stabilizable in the space $H^1(M)$
assuming that the Geometric Control Condition and 
some unique continuation condition are satisfied.

Recently, the authors proved in \cite{RZ2007b} that 
the cubic Schr\"odinger equation on the torus $\T$ with a localized
control 
\be
\label{cubic} 
iu_t + u_{xx} + \lambda |u|^2 u = ia(x)h(x,t),\quad  x\in \T, 
\ee
is locally exactly controllable in $H^s(\T )$ for all 
$s\ge 0$ (hence, in $L^2(\T )$). Inspired by the work of Russell-Zhang 
in \cite{rz-1}, the  method of proof combined the momentum approach and
Bourgain analysis. In the same paper, the local stabilization 
by the feedback law $h=a(x) u(x,t)$ was established by applying 
the contraction mapping theorem in some Bourgain space.
Finally, similar results were obtained with
Dirichlet (resp. Neumann) homogeneous boundary conditions thanks to 
an extension argument. 
More recently, Laurent has shown in \cite{laurent} that the system
(\ref{cubic}) is semiglobally exactly controllable and stabilizable. 
The same result has also been derived by Laurent in \cite{laurent2} for 
certain  manifolds of dimension 3, including $\T ^3$, $S^3$, 
and $S^2\times S^1$. 
The propagation of compactness and regularity proved in 
\cite{laurent,laurent2} plays a crucial role in the derivation of the 
stabilization results in these papers.
See also \cite{LRZ} for another application of these ideas to
the semiglobal stabilization of the periodic Korteweg-de Vries equation.

In addition, the authors considered in \cite{RZ2008} 
the following nonlinear Schr\"odinger equation
$$
iu_t + \Delta u + \lambda |u|^2 u =0
$$
posed on a bounded domain $\Omega $ in $\R ^n$ with either 
the Dirichlet boundary conditions or the Neumann boundary conditions.  
They showed that  if
\[ s>\frac{n}{2}, \]
or
\[ 0\leq s< \frac{n}{2} \mbox{ with } 1\leq n < 2+2s,\]
or 
\[ s=0,1 \mbox{ with } n=2,\]
then the systems  with control inputs acting on the whole boundary of
$\Omega$ are locally exactly controllable in the classical Sobolev 
space $H^s (\Omega) $ around any smooth solution of the  
Schr\"odinger equation.

The aim of this paper is to extend the results of \cite{RZ2007b} to
any dimension. More precisely, we shall assume that the spatial variable
lives in the rectangle 
$$ 
\Omega =(0,l_1)\times \cdots \times (0,l_n).
$$
We shall investigate the control properties of the semilinear 
Schr\"odinger equation 
\be
iu_t + \Delta u +\lambda |u|^{\alpha }u = i a(x)h(x,t),
\ee
where $\lambda \in \R$ and $\alpha\in 2\N^*$,
by combining new linear controllability results in the spaces
$H^s(\Omega )$ with Bourgain analysis. Let us briefly review
the results proved in this paper.


The internal controllability of the linear Schr\"odinger equation
on $\TN$ 
\be
\label{LS}
iu_t + \Delta u = i a(x)h(x,t), \quad x\in \TN,\  t\in (0,T)
\ee  
is established in $H^s(\TN)$ for any $s\ge 0$ and 
any function $a\not\equiv 0$. (Note that the Geometric Control Condition 
is not required.) It is derived from a well-known
result in $L^2(\TN)$, due to Jaffard \cite{jaffard} when $n=2$
and Komornik \cite{komornik} for any $n\ge 2$, by an argument allowing 
to shift the (state and control) space from $L^2(\TN )$ to $H^s(\TN )$. 
In particular,
the exact controllability in $H^s(\TN )$ will require a control input 
$h\in L^2(0,T;H^s(\TN))$.  Similar results with Dirichlet or Neumann
homogeneous boundary conditions are deduced by 
using the extension argument from \cite{RZ2007b}. 

The boundary controllability of the linear Schr\"odinger equation is considered
both with Dirichlet control 
\be
\label{dirichletN}
u=1_{\Gamma _0}h(x,t)
\ee
and with Neumann control
\be
\label{NeumannN}
\frac{\partial u}{\partial \nu}= 1_{\Gamma _0}h(x,t).
\ee
In \eqref{dirichletN} and in \eqref{NeumannN}, $\Gamma_0$ denotes an open
set in $\partial \Omega$. For the Dirichlet control, we shall prove that 
in {\em any} dimension $n\ge 2$ the exact controllability holds in 
$H^{-1}(\Omega )$ whenever  $\Gamma_0$ is a neighborhood of a vertex of
$\Omega$. The observability inequality for this 
(arbitrarily small) control region is actually 
derived from the corresponding observability
inequality for internal control by multiplier techniques. 

For the Neumann control, the exact controllability 
in $L^2(\Omega )$ is obtained
in any dimension when $\Gamma _0$ is a side. Finally, the 
results with Dirichlet (resp Neumann) boundary control
are extended to any Sobolev space $H^s(\Omega)$ with $s<1/2$ 
(resp. $s<1$)
by considering control inputs more regular in time, namely 
$h\in H^{\frac{s+1}{2}}(0,T;L^2(\partial \Omega ))$
(resp. $h\in H^{\frac{s}{2}}(0,T;L^2(\partial \Omega ))$). 

The extension of the above exact controllability results to the semilinear
Schr\"odinger equation 
\be
\label{NLS}
iu_t+\Delta u +\lambda |u|^{\alpha} u = i a(x)h(x,t)
\ee
is performed on the basis of Bourgain analysis. The needed linear
and multilinear estimates are combined with a fixed-point argument to produce
local exact controllability results. Sharp results (for the support
of the control input) are given for the internal control. Boundary
controllability results are derived from those established for the 
linear equation with the aid of estimates in Bourgain spaces of solutions
of boundary-value problems with boundary terms given by HUM.

Finally, the local exponential stabilization with an internal feedback law 
is proved by following the same approach as in \cite{RZ2007b}. 

The paper is organized as follows. The controllability  results
for the linear  Schr\"odinger equation are collected in Section 2. 
Section 3 is devoted to the controllability
of the semilinear equations. Section 4 deals with the 
internal stabilization issue.  
Multilinear estimates for nonlinearities
of the form $u^{\alpha _1}{\overline u}^{\alpha _2}$ are established in 
Appendix.

\section{Linear systems}

\subsection{Internal control}

We first consider the linear open loop control system  for the Schr\"odinger
equation posed on $\TN:=(-\pi ,\pi )^n$ with periodic boundary conditions:
\begin{equation}
\label{I1}
iu_t+\Delta u =iGh:=ia(x)h(x,t),\ \  u(x,0)=u_0(x),
\end{equation}
where $a\in C^\infty (\TN)$ is a given smooth real-valued function 
and $h=h(x,t)$ is the control input.

We denote by $H^s(\T ^n)$ the Sobolev space of the functions $u$ defined
on the torus $\TN$ (i.e. defined on $\R ^n$ and periodic of period
$2\pi$ with respect to each variable $x_i$)
for which the $H^s$ norm
$$
||u||_s=||(1-\Delta )^{s/2} u||_{L^2(\TN)}
$$
is finite.

We first establish an internal observability inequality for the solution 
$v(t)=W(t)v_0$ of
\begin{equation}
\label{A0}
\left\{
\begin{array}{l}
iv_t +\Delta v=0 \qquad (x,t)\in \TN\times \R ,\\
v(0)=v_0.
\end{array}
\right.
\end{equation}

\begin{prop}
\label{prop1}
{\bf (Observability inequality in $H^{-s}(\T ^n)$)}
Let $a\in C^\infty (\TN)$ with $a\ne 0$ and $T>0$.
Then for any $s\ge 0$ there exists a constant $c>0$ such that
for any  solution $v$ of
\eqref{A0} with $v_0\in H^{-s}(\T ^n)$, it holds
\begin{equation}
\label{A1}
||v_0||^2_{-s} \le c \int_0^T||av(t)||^2_{-s}dt.
\end{equation}
\end{prop}

\noindent{\em Proof.} We proceed in several steps.\\

\noindent
{\em Step 1.}  Assume that $s=0$, and let
$$
\omega =\{ x\in (-\pi ,\pi)^n;\ |a(x)|>||a||_{L^\infty(\TN)} /2 \}.
$$
Then, by \cite[Lemma 8.9]{KL},
there exists some  positive constant $c$ such that for any
square-summable sequence
$(c_k)_ {k\in \Z ^n \setminus \{ 0 \}}$ we have
\begin{equation}
\label{A2}
\sum_{k\ne 0}|c_k|^2
\le c\int_0^T\!\!\!\int_\omega \left|\sum_{k\ne 0}
c_k e^{ i (k\cdot x -|k|^2t)}\right|^2 dxdt.
\end{equation}
The result is still valid when the set of indices is changed into
$\Z ^n$  by \cite[Proposition 8.4]{KL}. This yields  \eqref{A1} when $s=0$. \\

\noindent
{\em Step 2.} We prove the weaker inequality
\begin{equation}
||v_0||^2_{-s} \le c
\left( \int_0^T || av(t)||^2_{-s}dt + ||v_0||^2_{-s-1}
\right)
\label{A20}
\end{equation}
by contradiction.
If \eqref{A20} is false, then there exists a sequence
$\{ v_j\}$ of solutions of
\eqref{A0} in $C([0,T];H^{-s}(\TN ))$ such that
\begin{equation}
1=||v_j(0)||^2_{-s} \ge j
\left( \int_0^T || av_j(t)||^2_{-s}dt + ||v_j(0)||^2_{-s-1}
\right) .
\label{A3}
\end{equation}
Since $v_j$ is bounded in $L^\infty ([0,T];H^{-s}(\TN ))$ and $(v_j)_t$ is
bounded in $L^\infty ([0,T];H^{-s-2}(\TN ))$ by \eqref{A0}, we infer
from Aubin's lemma that, for a subsequence again
denoted by $\{ v_j\}$, we have for $j\to \infty$
$$
\left\{
\begin{array}{ll}
v_j\to v \qquad &\text{ in } L^\infty ([0,T];H^{-s}(\TN )) 
\quad \text{weak }*\\
v_j\to v \qquad &\text{ in } C([0,T];H^r(\TN ))\quad \forall r<-s
\end{array}
\right.
$$
where $v\in C_w([0,T];H^{-s}(\TN ))$ is a solution of \eqref{A0}.
In particular, $v_j(0)\to v(0)$ in $H^r(\TN )$ for any $r<-s$. Since
$v_j(0)\to 0$ in  $H^{-s-1}(\TN )$ by \eqref{A3}, we conclude that $v\equiv 0$.
Let $w_j=(1-\Delta )^{-s/2}v_j$. Then $w_j\in L^\infty ([0,T];L^2(\TN))$ and
 $$
\left\{
\begin{array}{ll}
w_j\to 0 \qquad &\text{ in } L^\infty ([0,T];L^2(\TN))\quad \text{weak }* \\
w_j\to 0 \qquad &\text{ in } C([0,T];H^r (\TN ))\quad \forall r<0.
\end{array}
\right.
$$
Let us split $aw_j$ into
$$
aw_j=(1-\Delta )^{-s/2}(av_j) -(1-\Delta )^{-s/2}[a,(1-\Delta )^{s/2}]w_j.
$$
As the pseudodifferential operator $[a,(1-\Delta )^{s/2}]$ maps continuously
$H^r(\TN )$ into $H^{r-s+1}(\TN )$, we have that
\begin{equation}
\label{A4}
(1-\Delta )^{-s/2} [a,(1-\Delta )^{s/2} ]w_j \to 0  \qquad
\text{ in } C([0,T];H^r(\TN ))\ \text{ for any } r<1.
\end{equation}
Therefore, using \eqref{A3} and \eqref{A4}, we obtain that
$$
aw_j \to 0 \text{ in } L^2([0,T]; L^2(\TN)).
$$
Clearly, $w_j$ satisfies also the linear Schr\"odinger equation \eqref{A0},
so we infer from the observability inequality \eqref{A1}
established for $s=0$ that
$$
w_j(0)\to 0 \text{ in } L^2(\TN).
$$
It follows that $v_j(0)=(1-\Delta )^{s/2} w_j(0)\to 0$ in $H^{-s}(\TN )$,
contradicting the fact that $||v_j(0)||_{-s}=1$ for all $j$.\\

\noindent{\em Step 3.} We prove \eqref{A1} by contradiction.
If \eqref{A1} is false, there exists a sequence $\{ v_j\}$ of solutions
of \eqref{A0} in $C([0,T];H^{-s}(\TN ))$ such that
\begin{equation}
\label{A5}
1=||v_j(0)||^2_{-s}\ge j\int_0^T||av_j(t)||^2_{-s}dt \qquad
\forall j\ge 0.
\end{equation}
Extracting a subsequence if needed, we may assume that
\begin{eqnarray}
v_j\to v \qquad &&\text{ in } L^\infty ([0,T];H^{-s}(\TN )) \quad \text{weak}*
\label{A6}\\
v_j\to v \qquad &&\text{ in } C([0,T];H^r(\TN ))\quad \forall r<-s
\label{A7}
\end{eqnarray}
for some solution $v\in C_w([0,T];H^{-s}(\TN ))$ of \eqref{A0}, where
$C_w([0,T];H^{-s}(\T ^n)$ denotes the space of weakly sequentially 
continuous functions from $[0,T]$ to $H^{-s}(\T ^n)$ 
(see \cite[Lemme 8.1]{LM}).
Clearly, $a v_j\to av$ in $L^\infty ([0,T];H^{-s}(\TN ))$ weak $*$ which,
combined to \eqref{A5}, yields $a v\equiv 0$.
An application of Holmgren theorem (see e.g. \cite[Theorem 8.6.5]{hormander})
gives $v\equiv 0$. On the other hand,
\eqref{A7} gives $v_j(0)\to 0$ in $H^{-s-1}(\TN )$.
It then follows from \eqref{A20} that $v_j(0)\to 0$ in $H^{-s}(\TN )$,
and this contradicts \eqref{A5}.
\qed

Applying  HUM \cite{L} with $L^2(\T ^n)$ as pivot space, we infer from 
Proposition \ref{prop1} the following internal controllability 
of the linear Schr\"odinger equation in $H^s(\T ^n)$.
 
\begin{thm}
\label{thm1}
Let $T>0$ and $s\ge 0$ be given. Then for any 
$(u_0,u_1)\in H^s(\T ^n )\times H^s(\T ^n)$ there exists
a control $h\in L^2([0,T];H^s(\T ^n))$ such that the system  \eqref{I1}
admits a unique solution
$u\in C([0,T]; H^s (\T ^n))$ satisfying   $u(T)=u_1$.
Moreover, we can define a bounded operator
$$\Phi  :H^s(\T ^n)\times H^s(\T ^n)\to L^2([0,T];H^s(\T ^n))
$$
such that for any $(u_0,u_1)\in H^s(\T ^n)\times H^s(\T ^n)$ it holds
\begin{equation}
\label{oper}
W(T)u_0+\int_0^T W(T-\tau ) 
( G(\Phi (u_0,u_1)))(\cdot , \tau )d\tau =u_1.
\end{equation}
\end{thm}

\begin{figure}[hbtp]
\begin{center}
\epsfig{file=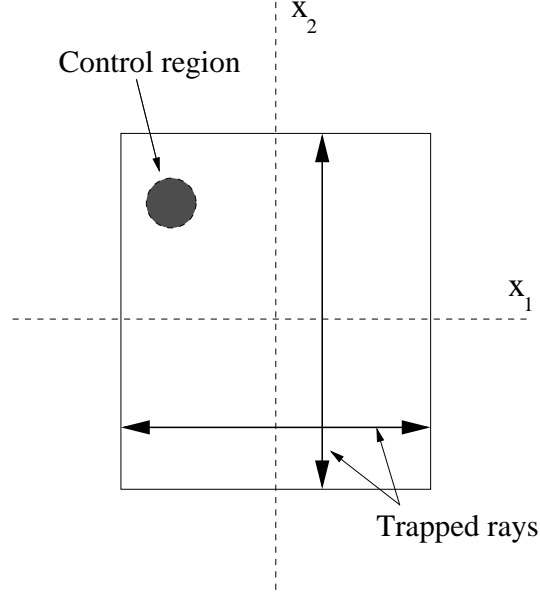,width=7cm}
\caption{Internal control of the Schr\"odinger equation.}
\label{tore1}
\end{center}
\end{figure}

The (small) control region is represented in Figure \ref{tore1}.
Trapped rays are drawn to mean that the wave equation fails to be 
controllable with such control regions.

\subsection{Boundary control}

In this section $\Omega =(0,\pi )^n$, and $\Gamma _0$ 
denotes an open set in $\partial \Omega$. 

\noindent
\subsubsection{Dirichlet boundary control}

We first adopt the following definition.

\begin{defi}
\label{def1}
The open set $\Gamma _0\subset\partial \Omega$ is called a {\em Dirichlet
control domain} if given any $u_0,\ u_1\in H^{-1}(\Omega)$ and
any time $T>0$, one may find a control $h\in L^2(0,T;L^2(\Gamma _0))$ 
such that the solution $u=u(x,t)$ of
\begin{equation}
\label{B0}
\left\{
\begin{array}{ll}
iu_t+\Delta u = 0\qquad &\text{ in }  \Omega\times (0,T)\\
u=1_{\Gamma _0}h(x,t)\qquad &\text{ on } \partial\Omega \times (0,T) \\
u(0)=u_0&
\end{array}
\right.
\end{equation}
satisfies $u(T)=u_1$.
\end{defi}

The following result provides Dirichlet control domains which are arbitrary
small in {\em any} dimension $n\ge 2$. Note that the wave equation 
fails to be controllable with such control domains.

\begin{figure}[hbtp]
\begin{center}
\epsfig{file=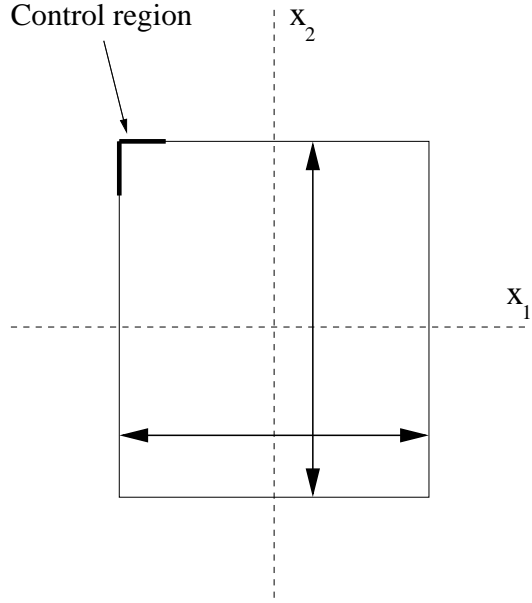, width=7cm}
\caption{Boundary control of the Schr\"odinger equation.}
\label{tore2}
\end{center}
\end{figure}

\begin{thm}
\label{thm2}
Let $\Omega =(0,\pi )^n$, and 
let $\Gamma _0\subset \partial \Omega$ be any open set
containing a vertex of $\partial \Omega$. Then $\Gamma _0$ is a Dirichlet
control domain.
\end{thm}

By Dolecki-Russell test of controllability (or HUM),
Theorem \ref{thm2} is a direct consequence of the following 
boundary observability result for the system
\begin{equation}
\label{B1}
\left\{
\begin{array}{ll}
iv_t+\Delta v = 0\qquad &\text{ in }  \Omega\times (0,T)\\
v=0\qquad &\text{ on } \partial\Omega \times (0,T) \\
v(0)=v_0.&
\end{array}
\right.
\end{equation}

\begin{prop}
\label{prop2}
Assume that the (open) control region $\Gamma _0\subset \partial \Omega$
contains a vertex of $\partial \Omega$. Then  for every $T>0$, there exists a
constant  $c>0$ such that
\be
\label{B7}
||\nabla v_0||^2_{L^2(\Omega )}
\le c\int_0^T\!\!\! \int_{\Gamma _0}
\left\vert\frac{\partial v}{\partial \nu} 
\right\vert ^2d\sigma dt
\ee
for any solution $v$ of \eqref{B1} with $v_0\in H^1_0(\Omega )$.
\end{prop}

\noindent{\em Proof.} We proceed in several steps.\\
{\em Step 1.} First, we prove an observability inequality in $H^1_0(\Omega)$
with an internal observation in an arbitrary subdomain of $\Omega$.
\begin{lem}
\label{lem1}
Let $\omega \subset \Omega$ be an arbitrary nonempty open set. 
Then there exists a constant $c>0$ such that
\be
\label{B6}
||\nabla v_0||^2_{L^2(\Omega )}
\le c\int_0^T\!\!\! \int_{\omega}
\left\vert\nabla v (x,t)\right\vert ^2 dxdt
\ee
for every solution $v$ of \eqref{B1} with 
$v_0\in H^1_0(\Omega )$.
\end{lem}
\noindent
{\em Proof of Lemma \ref{lem1}.} Extend $v$ to $(-\pi,\pi)^n\times (0,T)$ in
such a way that $v$ is an odd function of $x_i$ for each $i=1,...,n$, and
extend the initial state $v_0$ in a similar way. Then $v$ solves
\eqref{A0}. Writing $v_0=\sum_{k\in \Z ^n}c_ke^{ik\cdot x}$, we have that
$$
\nabla v(x,t)=\sum_{k\in \Z ^n} ic_k e^{i(k\cdot x-|k|^2t)}k.
$$
It follows then from \eqref{A2} that
\begin{eqnarray*}
||\nabla v_0||^2_{L^2(\TN )}
&=& \sum_{j=1}^n \sum_{k\in \Z ^n} |k_j|^2|c_k|^2  \\
&\le& c\sum_{j=1}^n \int_0^T\!\!\!\int_\omega
\big\vert \sum_{k\in \Z ^n} c_ke^{i(k\cdot x - |k|^2t)}k_j\big\vert ^2dxdt \\
&\le& c \int_0^T\!\!\!\int_\omega |\nabla v|^2 dxdt.
\end{eqnarray*}
The lemma is proved.\qed

\noindent
{\em Step 2.} We use the multiplier method to reduce the boundary
observation inequality to an internal observation inequality. Without
loss of generality, we may assume that $\Gamma_0$ is a (small)
neighborhood of the vertex $M=(\pi, ..., \pi)$ defined as
$$
\Gamma _0=\{ x\in \partial \Omega;\ 
x_1 + \cdots  + x_n > n \pi -\ep \},
$$
where $\ep$ is a (possibly small) positive number.
The following lemma is needed.
\begin{lem}
\label{lem2}
There exists a nonnegative function $\theta \in C^3(\R ^n)$
which is null on $\{ x\in\R ^n;\  x_1\le 0\}$ and strictly convex on
$(0,+\infty )^n\cap B_1(0)$.
\end{lem}
\noindent
{\em Proof of Lemma \ref{lem2}.}
Set $y^+=\max (y,0)$ for all $y\in \R$. Let
$$
\theta (x_1,...,x_n)=(x_1^+)^4
\big( 1+\delta \sum_{j=2}^n(x_j^+)^4 \big)
$$
where $\delta >0$ is a small number whose value will be specified later.
Clearly, $\theta$ is a nonnegative function 
of class $C^3$ on $\R ^n$, which vanishes on the set $\{ x_1\le 0\}$.
To prove that $\theta$ is strictly convex on $(0,+\infty)^n
\cap B_1(0)$, it is sufficient
to check that the Hessian matrix
\be
\label{hessian}
H(x)=\left( \frac{\partial ^2\theta }{\partial x_i\partial x_j} (x)\right)
\ee
is positive definite  for every $x\in (0,+\infty )^n\cap B_1(0)$.  
Simple computations give that for any $\xi\in\R ^n$,
$$
\xi ^TH(x)\xi =12 x_1^2 (1+\delta \sum_{j=2}^n x_j^4)\xi _1^2
+12\delta x_1^4 \sum_{j=2}^n x_j^2 \xi _j^2
+32\delta x_1^3\xi _1\sum_{j=2}^n x_j^3\xi _j.
$$
From Young inequality, we obtain that
$$
32 |x_1^3x_j^3\xi _1\xi _j|\le 26 x_1^2 x_j^4 \xi _1 ^2 
+10x_1^4 x_j^2 \xi _j^2,
$$
therefore
\be
\xi ^TH(x)\xi \ge  (12-26(n-1)\delta ) x_1^2 \xi_1 ^2
+2\delta  x_1^4 \sum_{j=2}^n x_j^2 \xi _j^2 \ge c|\xi |^2
\ee
if $x\in (0,+\infty)^n\cap B_1(0)$ and $\delta <(6/13)(n-1)^{-1}$. \qed
At this position, we need an identity from \cite{mach-2}.
\begin{lem}\cite[Lemma 2.2]{mach-2}
\label{lem3}
For any $q\in H^2(\Omega ,\R ^n)$ and any solution $v$ of
\eqref{B1} issued from $v_0\in H^1_0(\Omega )$, it holds
\begin{eqnarray}
&&\frac{1}{2}\int_0^T\!\!\!\int_{\partial \Omega}(q\cdot \nu)
\left\vert \frac{\partial v}{\partial \nu}\right\vert ^2 d\sigma dt 
= \frac{1}{2}\text{\rm Im } 
\int_{\Omega}(v q\cdot \nabla \bar v)dx\vert _0^T
\nonumber\\
&& +\frac{1}{2}\text{\rm Re }
\int_0^T\!\!\!\int_{\Omega}(v\nabla (\text{div }q)\cdot 
\nabla \bar v )dxdt
+\text{\rm Re } \int_0^T\!\!\!\int _{\Omega}\sum_{j,k=1}^n
\frac{\partial q_k   }{\partial x_j}\,
\frac{\partial \bar v}{\partial x_k}\,
\frac{\partial v     }{\partial x_j}\, dxdt.     \label{A10}
\end{eqnarray}
\end{lem}
Let
$$
\omega =\{ x\in \Omega;\ x_1+\cdots + x_n > n\pi -\ep \}.
$$
We readily infer from  Lemma \ref{lem2} that there exists a
convex function $\theta \in C^3(\overline{\Omega})$ which is
strictly convex on $\omega$ and null on $\overline{\Omega \setminus \omega}$.
Using \eqref{A10} with $q=\nabla \theta$ we obtain
\be
\label{A11}
\int_0^T\!\!\!\int_\omega
\nabla {\bar v}(x)^TH(x)\nabla v (x)\, dxdt
\le c\int_0^T\!\!\!\int _{\Gamma _0}  \left\vert
\frac{\partial v}{\partial \nu }\right\vert ^2d\sigma dt +
C_\delta \int_{\Omega} |v_0|^2 dx +\delta \int_\Omega
|\nabla v_0|^2 dx,
\ee
where $\delta >0$ is a small number and $H(x)$ denotes the Hessian matrix 
given in \eqref{hessian}.
In \eqref{A11}, we used the fact that both quantities
$||v(t)||_{L^2(\Omega )}$ and $||\nabla v(t)||_{L^2(\Omega )}$ are conserved.
Using Lemma \ref{lem1} and the fact that the Hessian matrix
$H(x)=(\partial ^2\theta/\partial x_i\partial x_j)(x)$ 
is positive definite on $\omega$, we obtain
\be
\label{A40}
||\nabla v_0||^2_{L^2(\Omega )}
\le c\int_0^T\!\!\!\int_{\Gamma _0}
\left\vert \frac{\partial v}{\partial \nu}\right\vert ^2 
d\sigma dt + C_\delta   \int_\Omega |v_0|^2dx.
\ee
for a convenient choice of $\delta$. The proof of the estimate
\be
\label{A41}
||v_0||^2_{L^2(\Omega)} \le c\int_0^T\!\!\! \int_{\Gamma _0}
\left\vert \frac{\partial v}{\partial \nu}\right\vert ^2 
d\sigma dt \ee
is classical (see e.g. \cite[pp. 27-28]{mach-2}). Then
\eqref{B7} follows from \eqref{A40}-\eqref{A41}. This completes
the proof of Proposition \ref{prop2} and of Theorem \ref{thm2}.
\qed

\begin{rem}
\begin{enumerate}
\item Theorem \ref{thm2} is stated for a square $\Omega =(0,\pi )^n$,
but it is valid (with the same proof) for any rectangle 
$\Omega =(0,l_1)\times \cdots \times (0,l_n)$. 
\item Using a frequential criterion and number theoretic arguments, 
Ramdani et al. \cite{RTTT} proved that when $n=2$,  
$\Gamma _0\subset \partial \Omega$ is a Dirichlet control 
domain if and only if $\Gamma _0$
has both a horizontal and a vertical components. It is however 
unclear whether the approach in \cite{RTTT} can yield 
a similar result for $n\ge 3$.
\item Using Theorem \ref{thm1} on a rectangle $\tilde \Omega
=(-1 ,\pi )\times (0,\pi )^{n-1}$ with a control input
supported in $\tilde \Omega \setminus \Omega$, and next taking the 
restriction  to $\Omega$, we infer that the linear Schr\"odinger 
equation is
controllable in $L^2(\Omega )$ with a control supported on a side.
(This fact can also be deduced from the Carleman inequalities
established in \cite{mor}.) This suggests that the condition for a domain
to be a Dirichlet control domain is less restrictive when the state
space is smoothed.  
\end{enumerate}
\end{rem}

We now aim to extend Theorem \ref{thm2} to a control result in a space
$H^s(\Omega)$, with $s\ge -1 $. We define 
$H^s_D(\Omega )=D(A_D^{\frac{s}{2}})$, where 
$A_D$ is the Dirichlet Laplacian; i.e., $A_Du=-\Delta u$ with domain
$D(A_D)=H^2(\Omega)\cap H^1_0(\Omega ) \subset L^2(\Omega )$.
We first need to replace the characteristic function $1_{\Gamma _0}$ by a 
smooth controller function $g\in L^\infty (\partial \Omega )$.  We adopt the following
\begin{defi}
Let $g\in L^\infty (\partial \Omega )$. We say that $g$ is a 
{\em smooth Dirichlet controller} if 
\begin{enumerate}
\item[(i)] there exists a constant $C>0$ such that 
\begin{equation}
\label{AA1}
||\nabla v_0||^2_{L^2(\Omega )} \le 
C\int_0^T\!\!\!\int_{\partial \Omega } g(x) 
\left\vert \frac{\partial v}{\partial \nu}\right\vert ^2 d\sigma dt 
\end{equation}
for any solution $v$ of \eqref{B1} emanating from $v_0\in H^1_0(\Omega )$
at $t=0$; 
\item for any face $F$ of $\partial \Omega $, $g_F=g_{\vert _F}\in C^\infty (F)$ 
and for all $k\ge 0$ 
\begin{equation}
\label{AA10}
\frac{\partial ^{2k+1} g_F}{\partial \nu ^{2k+1}} =0 \quad \text{ on } 
\partial F.
\end{equation}
\end{enumerate}
\end{defi}
Note that for any nonempty open set $\Gamma _0\subset \partial \Omega$ one
can construct a smooth Dirichlet controller $g$ supported in $\Gamma _0$. 
Consider for example a small neighborhood $\Gamma _0 = [0,\varepsilon ]^n\cap 
\partial \Omega $ of $0$ in $\partial \Omega$. A smooth Dirichlet controller
$g$ supported in $\Gamma _0$ is given by 
$$
g(x_1,...,x_n)=\prod_{i=1}^n \rho (x_i)
$$
where $\rho \in C^\infty (\R )$ fulfills
$$
\rho (s) = 
\left\{ 
\begin{array}{ll}
1 & \text{\rm if }  s\le \frac{\varepsilon}{4},\\
0 & \text{\rm if }  s\ge \frac{\varepsilon}{2}.
\end{array}
\right.
$$
Note also that $g\in C^0(\partial \Omega )$ and that the set 
$\{ x\in \partial \Omega;\  g(x)>0 \}$
is an open neighborhood of $0$ in $\partial \Omega$.

Let $g$ be a smooth Dirichlet controller, and let $S$ denote the bounded 
operator $H^1_0 (\Omega )\to  H^{-1}(\Omega )$ defined by $Sv_0=u(T)$, where
$u=u(x,t)$ solves 
\begin{equation}
\label{AA2}
\left\{
\begin{array}{ll}
iu_t+\Delta u = 0\qquad &\text{ in }  \Omega\times (0,T)\\
u=g(x)h(x,t)\qquad &\text{ on } \partial\Omega \times (0,T) \\
u(0)=0&
\end{array}
\right.
\end{equation}
with $h(x,t)=(\partial v/\partial \nu ) (x,t)$, $v=W_D(t)v_0$ denoting 
the solution of
\begin{equation}
\label{AA2bis}
\left\{
\begin{array}{ll}
iv_t+\Delta v = 0\qquad &\text{ in }  \Omega\times (0,T)\\
v=0\qquad &\text{ on } \partial\Omega \times (0,T) \\
v(0)=v_0.
\end{array}
\right.
\end{equation}
Applying HUM, we infer from the observability inequality \eqref{AA1}
that $S$ is invertible. We shall prove that a similar result holds in more 
regular spaces. 

\begin{thm}
\label{thm3}
Pick any number $s\in [-1,\frac{1}{2} )$. Then $S$ is an isomorphism from 
$H_D^{s+2}(\Omega )$ onto $H^s_D(\Omega )$. More precisely, for any 
$T>0$ and any $u_T \in H^s_D(\Omega )$, if we set 
$h(x,t)=(\partial v/\partial \nu ) (x,t)$ where $v$ denotes the solution of \eqref{AA2bis} with $v_0=S^{-1}u_T$, then $v_0\in H^{s+2}_D(\Omega )$, 
$h\in H^{\frac{s+1}{2}}(0,T;L^2(\partial \Omega ))$, and 
the solution $u$ of \eqref{AA2} satisfies  $u\in C([0,T]; H^s_D(\Omega ))$ and 
$u(T)=u_T$.
\end{thm}
\noindent{\em Proof.} 
{\em Step 1.} Let us first check that $S^{-1}$ is a bounded operator from
$H^s_D(\Omega )$ into $H^{s+2}_D(\Omega )$ for $s\in [-1,\frac{1}{2})$.  
The result is already known for $s=-1$. Assume first that 
$ -1 < s < 0 $, and pick any $u_T \in H^s_D(\Omega )$ decomposed as 
$$
u_T(x) =\sum_{p\in (\N ^*)^n} u_{T,p} \sin(p_1 x_1)\cdots \sin(p_nx_n), 
$$ 
with $\sum_{p\in (\N ^*)^n} |p|^{2s}|u_{T,p}|^2<\infty$. 
Let $v_0 = S^{-1}(u_T)\in H^1_D(\Omega )$
decomposed as
\begin{equation}
\label{XY1}
v_0(x) =\sum_{p\in (\N ^*)^n} 
v_p \sin(p_1x_1)\cdots \sin(p_nx_n), 
\end{equation}
and let $v$ denote the solution of \eqref{AA2bis}. The control 
given by HUM driving \eqref{AA2} from $0$ to $u_T$ reads
\be
\label{h1992}
h(x,t):=\partial v/\partial \nu
=\sum_{p\in (\N ^*)^n} 
v_p e^{-i|p|^2t}
\frac{\partial}{\partial \nu}
(\sin(p_1x_1)\cdots \sin(p_nx_n)).
\ee
Let us write the solution $u=u(x,t)$ of \eqref{AA2} in the form
\begin{equation}
\label{XY2}
u(x,t) =\sum_{p\in (\N ^*)^n} u_p(t)
\sin(p_1x_1)\cdots \sin(p_nx_n). 
\end{equation}
The moments $\{ u_p(t)\} _{p\in (\N ^*)^n}$ can be computed from 
the control input $h$ by using duality. Scaling in \eqref{AA2} 
by $\overline{w}$, where $w=W_D(t)w_0$ is a smooth solution, we obtain
$$
i\int_\Omega u(x,t)\overline{w(x,t)}\, dx =
\int_0^t \int_{\partial \Omega} g(x) h(x,\tilde t) 
\overline{\frac{\partial w}{\partial \nu}}\, d\sigma (x) d\tilde t.  
$$ 
Pick any $q\in (\N ^*)^n$ and  choose
$w_0(x)=\sin ( q_1x_1)\cdots \sin(q_nx_n)$.
We obtain from \eqref{h1992} that 
\begin{eqnarray}
(\frac{\pi}{2})^n ie^{i|q|^2t} u_q(t)
&=& \int_0^t \int_{\partial \Omega} g(x)
h (x,\tilde t) e^{i|q|^2\tilde t }
\frac{\partial}{\partial\nu}(\sin (q_1 x_1)\cdots \sin (q_n x_n))
d\sigma (x)d\tilde t \nonumber \\
&=& \sum_{p \in (\N ^*)^n} v_p
(\int_0^t e^{i(|q|^2-|p|^2)\tilde t } d\tilde t) \nonumber \\
&&\!\!\!\!\!\! \times\!\!\!
\int_{\partial \Omega} g(x)
\frac{\partial}{\partial\nu}(\sin (p_1 x_1)\cdots \sin (p_n x_n))
\frac{\partial}{\partial\nu}(\sin (q_1 x_1)\cdots \sin (q_n x_n)) d\sigma (x).
\label{XY3}
\end{eqnarray}
It follows that for $t=T$
\be
\label{solution}
S(v_0)=u_T=u(T)=\sum_{q\in (\N ^*)^n} 
\left( \sum_{p\in (\N ^*)^n} a_{q,p} v_p \right) 
\sin (q_1 x_1)\cdots \sin (q_n x_n)
\ee
with 
\be
\label{AA3}
a_{q,p}=
-(\frac{2}{\pi}) ^{n} \frac{e^{-i|p|^2T} - e^{-i|q|^2T}}{|q|^2-|p|^2} 
\int_{\partial \Omega} g(x)
\frac{\partial}{\partial\nu}(\sin (p_1 x_1)\cdots \sin (p_n x_n)) 
\frac{\partial}{\partial\nu}(\sin (q_1 x_1)\cdots \sin (q_n x_n)) d\sigma (x).
\ee
In \eqref{AA3}, we used the convention that 
\begin{equation}
\label{XY4}
\frac{e^{-i|p|^2t} - e^{-i|q|^2t}}{|q|^2-|p|^2}
=it e^{-i|q|^2t} \qquad \text{ for }\ |p|=|q|.
\end{equation}
Introduce the operator $D^\sigma$  defined by 
$$
D^\sigma 
\left( 
\sum_{p\in (\N ^*)^n} c_p\sin (p_1x_1)\cdots \sin (p_nx_n) 
\right)
=\sum_{ p \in (\N ^*)^n} |p|^\sigma c_p\sin (p_1x_1)\cdots \sin (p_nx_n).
$$
In what follows, $\sum _p$ and $\sum _q$ will stand for  
$\sum _{p\in (\N ^*)^n}$ and   $\sum _{q\in (\N ^*)^n}$, respectively.
We aim to prove that $v_0\in H^{s+2}_D (\Omega )$ for 
$u_T \in H^s_D ( \Omega )$. For $v_0$ given by \eqref{XY1}, let
$$
||v_0||_s^2=\sum_p |p|^{2s} |v_p|^2.
$$
$C$ denoting a constant varying from line to line, we have that
\begin{eqnarray}
||v_0||_{s+2} 
&\le& ||D^{s+1} v_0||_1 \nonumber\\
&\le & C||S(D^{s+1}v_0)||_{-1} \nonumber\\
&\le & C \left( ||D^{s+1}(Sv_0)||_{-1}  + ||[S,D^{s+1}]v_0||_{-1}\right) 
\nonumber\\
&\le & C\left( ||u_T||_s + || [S,D^{s+1} ] v_0 ||_{-1} \right) . 
\label{ABCD1}
\end{eqnarray}
Clearly
$$
[S,D^{s+1}] v_0 = \sum_q \left( \sum _p a_{q,p}
(|p|^{s+1} - |q|^{s+1}) v_p\right) \sin (q_1x_1)\cdots \sin (q_nx_n),
$$
hence 
$$
||[S,D^{s+1}] v_0 ||_{-1}^2 = 
\sum _q |q|^{-2} [ \sum_p a_{q,p} (|p|^{s+1} - |q|^{s+1}) v_p ] ^2.
$$
Writing $\partial \Omega =\cup_{ 0\le l < 2^n -1 } F_l$, 
where the $F_l$'s denote the faces of $\Omega$, the integral term 
in \eqref{AA3} may be written $\sum_{0\le l < 2^n - 1 } I_{F_l}$, with
$$
I_{F_l}:= 
\int_{F_l} g(x)
\frac{\partial}{\partial\nu}(\sin (p_1 x_1)\cdots \sin (p_n x_n)) 
\frac{\partial}{\partial\nu}(\sin (q_1 x_1)\cdots \sin (q_n x_n)) d\sigma (x).
$$
Let us estimate $I_{F_l}$ for $F_0:=\{ x\in \partial \Omega ;  x_n=0\} = 
[0,\pi ]^{n-1}\times \{ 0 \} $. Then 
\begin{eqnarray*}
|I_{F_0}| &=& p_n q_n 
\left\vert 
\int_{[0,\pi ]^{n-1}} g(x_1,...,x_{n-1},0)
[\prod_{j=1}^{n-1}\sin (p_jx_j)\sin (q_jx_j)]
dx_1\cdots dx_{n-1}
\right\vert\\
&=& p_n q_n 
\left\vert 
\int_{[0,\pi ]^{n-1}} g(x_1,...,x_{n-1},0)
[\prod _{j=1}^{n-1}\frac{1}{2} 
(\cos (p_j-q_j)x_j - \cos (p_j+q_j)x_j ) 
] dx_1\cdots dx_{n-1} \right\vert.
\end{eqnarray*}
Using \eqref{AA10} and integrations by parts, we see that for every 
$k\in \N$, we have for some constant $C_k>0$
\begin{equation}
\label{XYZ}
|I_{F_0} | \le C_k p_nq_n \prod _{j=1}^{n-1} \langle p_j-q_j\rangle ^{-k}.
\end{equation}
The corresponding contribution in $||[S, D^{s+1}]v_0||_{-1}^2$
is therefore estimated by 
$$A_{F_0}=\sum _q |q|^{-2}
\left( 
\sum_p p_n q_n (\prod _{j=1}^{n-1} \langle p_j- q_j\rangle ^{-k}) 
\langle |q|^2 - |p|^2 \rangle ^{-1} ||p|^{s+1}-|q|^{s+1}| |v_p| 
\right) ^2.$$
Since 
$$
\frac{\vert |p|^{s+1} -|q|^{s+1}\vert }{\langle |q|^2 - |p|^2 \rangle}
\le C \frac{\left\vert |p| - |q| \right\vert (|p|^s +|q|^s)}
{\langle |q|^2 -|p|^2 \rangle}
\le C \frac{|p|^s + |q|^s}{ |p| + |q| }
$$
we have by Cauchy-Schwarz
\begin{eqnarray}
A_{F_0} &\le& 
C\sum_q 
[\sum _p p_n (\prod _{j=1}^{n-1} \langle p_j-q_j\rangle ^{-k})
\frac{|p|^s + |q|^s}{|p|+ |q|} |v_p| ]^2\nonumber \\
&\le& C\sum_q 
\left( \sum_p \frac{|p|^{2s} + |q|^{2s}}{(|p|+ |q|)^2}\prod_{j=1}^{n-1}
\langle p_j-q_j\rangle ^{-k} \right) \cdot 
\left( \sum _p p_n^2 |v_p|^2 \prod _{j=1}^{n-1}
\langle p_j - q_j \rangle ^{-k}  \right) \label{AA20}
\end{eqnarray}
Pick any $k>1$. Then, as $s<0$, if we choose $k>1$
$$
\sum_{q_n}\sum_p \frac{|p|^{2s} + |q|^{2s}}{(|p|+ |q|)^2}
\prod_{j=1}^{n-1} \langle p_j-q_j\rangle ^{-k} 
\le 
\sum_{q_n}\sum_{p_n} \frac{p_n^{2s} + q_n^{2s}}{(p_n+q_n)^2}
\sum_{p_1,\ldots ,p_{n-1}} \prod_{j=1}^{n-1}\langle p_j-q_j\rangle ^{-k}<\infty. 
$$
Therefore
\begin{eqnarray*}
 A_{F_0} &\le& C\sum_{q_1,\ldots ,q_{n-1}} \sum_p p_n^2 |v_p|^2 
\prod_{j=1}^{n-1} \langle p_j-q_j\rangle ^{-k} \\
&\le & C\sum_p |p|^2 |v_p|^2 \sum_{q_1,\ldots ,q_{n-1}}
\prod_{j=1}^{n-1} \langle p_j-q_j\rangle ^{-k} \\
&\le & C\sum_{p} |p|^2 |v_p|^2.   
\end{eqnarray*}
The estimate for another face $F_l$ is similar. We conclude that 
$$
||[S,D^{s+1}]v_0||^2_{-1} \le C ||v_0||_1^2 \le C ||u_T||_{-1}^2
$$ 
hence, with \eqref{ABCD1}, $v_0\in H_D ^{s+2}(\Omega )$. 
Let us now assume that $u_T \in H_D ^s (\Omega )$ with $0\le s<\frac{1}{2}$. 
The proof is carried out as above when $-1<s<0$, except for the estimate of 
$A_{F_0}$ in \eqref{AA20}. We know from the lines above that 
$v_0\in H_D ^\sigma (\Omega )$ for any $\sigma <2$. Then, by Cauchy-Schwarz
inequality, 
\begin{eqnarray}
A_{F_0} 
&\le& C\sum_q \left(\sum _p p_n \prod_{j=1}^{n-1}
\langle p_j-q_j\rangle ^{-k} \frac{|p|^s+|q|^s}{|p| + |q|} |v_p|  \right)^2
\label{orion1}\\
&\le& C\sum_q 
\left(\sum _p \frac{|p|^{2s}+|q|^{2s}}{(|p| + |q|)^2} 
|p|^{-1} \prod_{j=1}^{n-1}\langle p_j-q_j\rangle ^{-k}  \right)
\left(\sum _p p_n^2 |p||v_p|^2
\prod_{j=1}^{n-1}
\langle p_j-q_j\rangle ^{-k} \right) . \nonumber
\end{eqnarray}
Note that 
$$
\sum_{q_n}\left( \sum_p \frac{|p|^{2s}+|q|^{2s}}{(|p|+|q|)^2} |p|^{-1}
\prod _{j=1}^{n-1}\langle p_j-q_j\rangle ^{-k}\right)
\le C(S_1+S_2+S_3) 
$$
where
\begin{eqnarray*}
S_1&=& \sum_{q_n}\left( \sum_p 
\frac{ |p|^{2s-1} }{ (|p|+|q|)^2} \prod_{j=1}^{n-1} 
\langle p_j-q_j\rangle  ^{-k}\right) \\
S_2&=& \sum_{q_n}\left( \sum_p 
\frac{ q_n^{2s}|p|^{-1} }{ (|p|+|q|)^2} \prod_{j=1}^{n-1} 
\langle p_j-q_j\rangle  ^{-k}\right) \\
S_3&=& \sum_{q_n}\left( \sum_p 
\frac{ |q'|^{2s}|p|^{-1} }{ (|p|+|q|)^2} \prod_{j=1}^{n-1} 
\langle p_j-q_j\rangle  ^{-k}\right) \qquad \text{ where }\ q=(q',q_n).
\end{eqnarray*}
Since $2s-1<0$, 
$$
S_1\le \sum_{q_n}
\left( 
\sum_p \frac{p_n^{2s-1}}{(p_n+q_n)^2}
\prod_{j=1}^{n-1}\langle p_j-q_j\rangle ^{-k} 
\right) \le const<\infty .
$$
Also, 
\begin{eqnarray*}
S_2 &\le & \sum_{q_n}\left( 
\sum_p \frac{q_n^{2s}p_n^{-1}}{(p_n+q_n)^2} \prod_{j=1}^{n-1}
\langle p_j-q_j\rangle ^{-k}  \right) \\
&\le& C\sum_{p_n}\sum_{q_n} \frac{q_n^{2s} p_n^{-1} }{ (p_n+q_n)^2} \\
&\le& C\sum_{p_n}\left( \frac{1}{p_n(p_n+1)} + p_n^{2s-3} + 
\int_1^\infty 
\frac{x^{2s}}{p_n(p_n+x)^2} dx \right) \\
&\le& 
C\left( 1+\sum_{p_n\ge 1} p_n^{2s-2} \int_0^{+\infty} \frac{y^{2s}}{(1+y)^2}dy \right)\\
&\le& const<\infty.
\end{eqnarray*}
Finally, 
$$
S_3 \le |q'|^{2s} \sum_{q_n}\sum_p \frac{p_n^{-1}}{(p_n+q_n)^2}
\prod_{j=1}^{n-1} \langle p_j-q_j\rangle ^{-k} 
\le C |q'| ^{2s}.
$$
It follows that 
$$
A_{F_0} \le C \sum_{q_1,\ldots ,q_{n-1}} \sum_p p_n^2 |p| |v_p|^2 
|q'| ^{2s} \prod _{j=1}^{n-1} \langle p_j-q_j\rangle ^{-k}. 
$$
Note that 
$$
\sum_{q_1,\ldots ,q_{n-1}} |q'| ^{2s} 
\prod_{j=1}^{n-1} \langle p_j - q_j \rangle ^{-k} 
\le C |p'| ^{2s}  
$$
since, for $k>2s+1$, 
$$
\sum_{q_j} q_j^{2s} \langle p_j-q_j \rangle ^{-k} 
\le Cp_j^{2s}. 
$$
(Split the sum into one for $q_j\le 2p_j$, and another one for $q_j>2p_j$.)
Therefore, since $0\le s<1/2$, 
\begin{equation}
\label{orion2}
A_{F_0} \le C \sum_p |p|^{3+2s} |v_p|^2 =
||v_0||^2_{s+\frac{3}{2}} \le C ||u_T||^2_{s-\frac{1}{2}}.
\end{equation}
Thus, we have proved that $S^{-1}$ is bounded from $H^s_D(\Omega )$ into $H^{s+2}_D(\Omega )$ for $-1\le s<\frac{1}{2}$. Note that, for $v_0\in 
H^{s+2}_D(\Omega )$,
$h\in H^{\frac{s+1}{2}}(\T ; L^2(\partial \Omega ))$ 
by  \eqref{h1992}. \\
{\em Step 2.} Since $S$ is an isomorphism from $H^{1}_D(\Omega )$ onto
$H^{-1}_D(\Omega )$, it remains to prove that $S$ maps $H^{s+2}_D(\Omega )$ into
$H^s_D(\Omega )$. The proof of Theorem \ref{thm3} 
will thus be complete with the following result.
\begin{prop}
\label{regularite}
Let $s\in [-1,\frac{1}{2})$ and $T>0$. For any $v_0\in H^{s+2}_D(\Omega )$, 
let $u=\Gamma v_0$ denote the solution of  \eqref{AA2} associated with 
$h=\partial v/\partial \nu$, where $v(t)=W_D(t) v_0$. Then $\Gamma $ is a bounded
operator from $H^{s+2}_D(\Omega )$ into $C([0,T];H^s_D(\Omega ))$.   
\end{prop}
\noindent
{\em Proof of Proposition \ref{regularite}.}
It is well known that for any $h\in L^2(0,T;L^2(\partial \Omega ))$, 
there exists a unique solution $u\in C([0,T];H^{-1}(\Omega ))$ 
in the transposition sense of \eqref{AA2} (see e.g. \cite{mach-2}). 
The result is therefore true for $s=-1$.  Let us now assume that $s\in (-1,1/2)$. 
From Step 1, we know that $u$ is given by  
\begin{equation}
\label{XY5}
u(t)
= -\left( \frac{2}{\pi}\right) ^n 
\sum_{q\in (\N ^*)^n}
\left(
\sum_{p\in (\N ^*)^n} 
v_p\frac{ e^{-i|p|^2 t} - e^{-i|q|^2 t}}{|q|^2-|p|^2} 
I(g,p,q)\right) 
\sin(q_1x_1)\cdots \sin (q_nx_n)
\end{equation}
where 
\begin{equation}
\label{XY6}
I(g,p,q)=\int_{\partial \Omega} g(x)\frac{\partial}{\partial \nu}
(\sin (p_1x_1)\cdots \sin (p_nx_n))\frac{\partial}{\partial \nu}
(\sin (q_1x_1)\cdots \sin (q_nx_n)) d\sigma (x).
\end{equation}
Again  $I(g,p,q)=\sum_{0\le l<2^n-1}I_{F_l}$, where the $F_l$'s denote the faces
of $\Omega$ and $I_{F_l}$ is given in \eqref{XY3}.
We have that 
\begin{eqnarray*}
||\Gamma v_0||_{L^\infty (0,T;H^s_D (\Omega ))}
&= & ||D^{s+1}(\Gamma v_0)||_{L^\infty (0,T;H^{-1}_D(\Omega ))} \\
&\le& 
||\Gamma (D^{s+1}v_0)||_{L^\infty (0,T;H^{-1}_D(\Omega ))}
+||[\Gamma ,D^{s+1}]v_0||_{L^\infty (0,T;H^{-1}_D(\Omega ))}
\cdot
\end{eqnarray*} 
Since
$$
||\Gamma (D^{s+1} v_0)||_{L^\infty (0,T;H^{-1}_D(\Omega ))}
\le C ||D^{s+1}v_0||_1 \le C ||v_0||_{s+2},
$$
it remains to estimate the commutator $[\Gamma , D^{s+1}]v_0$.
Clearly
\begin{equation}
([\Gamma ,D^{s+1}]v_0)(t)
= -\left( \frac{2}{\pi}\right) ^n 
\sum_{q} \left(
\sum_{p;|p|\ne |q|} v_p \frac{|p|^{s+1}-|q|^{s+1}}{|q|^2-|p|^2}
(e^{-i|p|^2 t} - e^{-i|q|^2 t})
I(g,p,q)\right) 
\prod_{j=1}^n\sin(q_jx_j).
\end{equation}
The contribution in $||([\Gamma , D^{s+1}]v_0)(t)||_{-1}^2$
due to $F_0=\{ x\in\partial \Omega; \ x_n=0\}$ is estimated with 
\eqref{XYZ} by 
\begin{eqnarray*}
B_{F_0} &\le& C\sum_q |q|^{-2}\left( \sum_{p,|p|\ne |q|}
|v_p| \frac{|p|^s+|q|^s}{|p|+|q|}|I_{F_0}|\right) ^2\\
&\le& C\sum_q 
\left(\sum_{p;|p|\ne |q|} |v_p| \frac{|p|^s + |q|^s}{|p|+|q|} p_n \prod_{j=1}^{n-1}
\langle p_j-q_j\rangle ^{-k}  \right) ^2. 
\end{eqnarray*}
Therefore, using the estimation of the r.h.s. of \eqref{orion1} in 
\eqref{orion2}, we conclude that for $s<1/2$
$$
B_{F_0} \le C||v_0||^2_{s+\frac{3}{2}},
$$
the constant $C$ being uniform in $t\in [0,T]$. Therefore
$$
||[\Gamma, D^{s+1}]v_0||_{L^\infty (0,T;H^{-1}_D(\Omega ))} 
\le C ||v_0||_{s+2}\cdot
$$
Thus, we have proved that 
\begin{equation}
\label{PQR}
||u||_{L^\infty (0,T;H^s_D(\Omega ))} 
\le C ||v_0||_{H^{s+2}_D(\Omega )}\cdot
\end{equation}
Since $u\in C([0,T];H^{-1}_D(\Omega ))$, we conclude that 
$u\in C_w ([0,T];H^s_D(\Omega ))$. If we pick $\tilde s\in (s,1/2)$ 
and $\tilde v_0\in H^{\tilde s +2}_D(\Omega )$, the corresponding solution
$\tilde u$ belongs to $C_w ([0,T];H_D^{\tilde s}(\Omega ))$, hence to 
$C([0,T];H^s_D(\Omega ))$, the embedding
$H^{\tilde s}_D (\Omega ) \subset H^s_D (\Omega )$ being compact.
It follows from \eqref{PQR}
combined to the density of $H^{\tilde s +2}_D(\Omega )$ in 
$H^{s+2}_D(\Omega )$ that $u\in C([0,T];H^s_D(\Omega ))$ for
$v_0\in H^{s+2}_D(\Omega )$.  
In particular, $u(T)\in H^s_D(\Omega )$, so that $S$ is an isomorphism from 
$H^{s+2}_D(\Omega ) $ onto $H^s_D(\Omega )$. This completes the proof of 
Proposition \ref{regularite} and of Theorem \ref{thm3}.
\qed

\noindent
\subsubsection{Neumann boundary control}

We adopt the following definition.

\begin{defi}
The open set $\Gamma _0\subset\partial \Omega$ is called a {\em Neumann
control domain} if given any $u_0,\ u_1\in L^2(\Omega)$ and
any time $T>0$, one may find a control $h\in L^2(0,T;L^2(\Gamma _0))$ 
such that the solution $u=u(x,t)$ of
\begin{equation}
\label{B123}
\left\{
\begin{array}{ll}
iu_t+\Delta u = 0\qquad &\text{ in }  \Omega\times (0,T)\\
\frac{\partial u}{\partial \nu}=1_{\Gamma _0}h(x,t)\qquad &\text{ on } 
\partial\Omega \times (0,T) \\
u(0)=u_0&
\end{array}
\right.
\end{equation}
satisfies $u(T)=u_1$.
\end{defi}

The following result provides Neumann control domains 
in {\em any} dimension $n\ge 2$. 

\begin{prop}
\label{prop20}
Let $\Omega =(0,\pi )^n$, and 
let $\Gamma _0\subset \partial \Omega$ be a side of $\Omega$. 
Then $\Gamma _0$ is a Neumann control domain.
\end{prop}
\noindent
{\em Proof.} Assume e.g. that $\Gamma _0= \{ 0 \} \times (0,\pi )^{n-1}$. 
By Dolecki-Russell criterion, we only have to check the
following observability inequality 
\be
||v_0||^2_{L^2(\Omega )} \le C \int_0^T\int_{\Gamma _0}|v(x,t)|^2 d\sigma dt
\ee
where $v_0$ is any function in $L^2(\Omega )$ and $v=v(x,t)$ solves
\begin{equation}
\label{B124}
\left\{
\begin{array}{ll}
iv_t+\Delta v = 0\qquad &\text{ in }  \Omega\times (0,T)\\
\frac{\partial v}{\partial \nu}=0\qquad &\text{ on } 
\partial\Omega \times (0,T) \\
v(0)=v_0.&
\end{array}
\right.
\end{equation}
Expanding $v_0$ as
$$
v_0(x)=\sum_{k\in \N ^n} c_k\cos(k_1x_1)\cdots \cos (k_nx_n),
$$ 
then the corresponding solution $v(x,t)$ reads
$$
v(x,t)=\sum_{k\in \N ^n} c_ke^{-i|k|^2t}
\cos(k_1x_1)\cdots \cos (k_nx_n).
$$
It follows that
\begin{eqnarray*}
\int_0^T\int_{\Gamma _0}|v(x,t)|^2 d\sigma dt
&= &\int_0^T \int_{(0,\pi )^{n-1}} | \sum_{k\in \N ^n} c_ke^{-i|k|^2t}
\cos(k_2x_2)\cdots \cos (k_n x_n)|^2 dx_2\cdots dx_n dt \\
&\sim& \sum_{k_2,...,k_n\ge 0}
\int_0^T \left\vert \sum_{k_1\ge 0} c_ke^{-ik_1^2t}\right\vert ^2dt
\sim \sum_{k\in \N ^n} |c_k|^2 \sim ||v_0||^2_{L^2(\Omega )},
\end{eqnarray*}
where we used the orthogonality of the functions 
$\cos(k_2x_2)\cdots \cos(k_n x_n)$
in $L^2(\Gamma _0)$ and Ingham's lemma.\qed


We now aim to extend Proposition \ref{prop20} to a control result in a space
$H^s(\Omega)$, $s>0$. We define $H^s_N(\Omega )=D(A_N^{\frac{s}{2}})$, where 
$A_N$ is the Neumann Laplacian (i.e. $A_N u=u-\Delta u$ with 
$D(A_N)=\{ u\in H^2(\Omega),\ \partial u/\partial \nu = 0\text{ on } 
\partial \Omega\} \subset L^2(\Omega )$).
A result similar to Theorem \ref{thm3} may be obtained along the same lines.
We limit ourselves to giving a weaker result with a very short proof.

\begin{thm}
\label{thm4}
Let $\Gamma _0$ be a Neumann control domain, $T=2\pi$, $s \in [0,1)$ 
and $u_0,u_1 \in H^s_N(\Omega )$. 
Then there exists a control input 
$h\in H^{\frac{s}{2}}(\T ;L^2(\partial\Omega ))$ such that the solution 
$u$ of \eqref{B123} satisfies $u(T)=u_1$. 
\end{thm}
\noindent{\em Proof.}  Without loss of generality, we may assume that
$u_0=0$. 
A direct computation shows that for any (smooth) solution $u$ of 
\eqref{B123} emanating from $u_0=0$ and any (smooth) solution 
$v$ of \eqref{B124}, it holds
\be
\label{identity}
i\int_\Omega u(x,T)\overline{v(x,T)}\, dx
= - \int_0^T\!\!\!\int_{\partial \Omega} 1_{\Gamma_0}h(x,t)
\overline{v} d\sigma dt. 
\ee
As usual, for any $h\in L^2(0,T; L^2(\partial \Omega ))$, the solution 
$u\in C([0,T];L^2(\Omega ))$  of \eqref{B123} is defined by 
\be
i(u(t) , v(t) )_{L^2(\Omega )}
= - (h, 1_{\Gamma _0}v)
_{L^2(0,t;L^2(\partial \Omega))}, \qquad\forall t\in [0,T],\ \forall v_0\in 
L^2(\Omega)
\ee 
where $v(t)$ solves \eqref{B124}.\\
{\sc Claim 1.}  If $v_0\in H^{-s}_N(\Omega )$ for some $s\in \R$, then 
$v\in H^{-\frac{s}{2}}(\T ;L^2(\partial \Omega ))$. \\
Indeed, if we write
$v_0=\sum_{k\in \N ^n} c_k\cos (k_1x_1)\cdots\cos (k_n x_n)$ and
$$v(x,t)=\sum_{k\in \N ^n} c_k e^{-i|k|^2t}
\cos (k_1x_1)\cdots\cos (k_n x_n)
$$
then we have that 
\begin{equation}
\label{LQR}
||v||^2_{H^{-\frac{s}{2}}(\T ,L^2(\partial \Omega ))} 
\sim \sum_{k}(1+|k|^2)^{-s} |c_k|^2\sim ||v_0||^2_{H^{-s}_N(\Omega )}.
\end{equation} 
We may rewrite \eqref{identity} in the form
\be
\label{identitybis}
i\langle u(T) , v(T) \rangle _{H^s_N, H^{-s}_N}
= - \langle h, 1_{\Gamma _0} v
\rangle_{H^{\frac{s}{2}}(\T ;L^2(\partial \Omega)),
 H^{-\frac{s}{2}}(\T ; L^2(\partial \Omega))}.
\ee 
Note that $u\in C([0,T]; H^s_N(\Omega ))$ if $0\le s<1$. 
It remains to establish the following\\
{\sc Claim 2.} (Observability inequality) The following 
estimate holds for the solutions of \eqref{B124}:
\be
\label{obs}
||1_{\Gamma _0} v||^2
_{H^{-\frac{s}{2}}(\T ;L^2(\partial\Omega ))} \ge const||v_0||^2_{H^{-s}_N
(\Omega )  }
\ee
If \eqref{obs} is not true, one can construct a sequence $\{ v_j\}$ such that
\be
\label{obs2}
j||1_{\Gamma _0} v_j||^2
_{H^{-\frac{s}{2}}(\T ;L^2(\partial\Omega ))} < 
||v_j(0)||^2_{H^{-s}_N(\Omega ) }=1.
\ee
Let $w_j=(1-\partial _t^2)_p^{-\frac{s}{4}} v_j$, where for 
any $\sigma\in \R$
$$
(1-\partial _t^2)_p^{\sigma} \sum_{l\in \Z}c_l e^{ilt} =
\sum_{l\in \Z}(1+|l|^2)^\sigma c_l e^{ilt}.
$$
 Then $w_j$ solves
\eqref{B124} with $w_j(0)$ substituted to $v_0$, and from \eqref{obs2} we obtain
\be
\label{obs3}
1_{\Gamma _0} w_j \to 0 \ \text{ in } \ 
L^2(\T ;L^2(\partial \Omega )).
\ee
As $\Gamma _0$ is a Neumann control domain, we infer that 
$w_j(0)\to 0$ in $L^2(\Omega)$, hence 
$$
w_j\to 0 \text{ in } L^2(\T ; L^2 (\partial \Omega )).
$$
This gives
$$
v_j\to 0 \text{ in } H^{-\frac{s}{2}}(\T ;L^2(\partial \Omega)).
$$
Using \eqref{LQR}, we infer that 
$v_j(0) \to 0$ in $H^{-s}_N(\Omega )$, which contradicts
\eqref{obs2}. 
This completes the proof of Theorem \ref{thm4}.
\qed

\section{Nonlinear systems}

\subsection{Internal control}
In this section we consider the following
nonlinear control system
\begin{equation}
\left \{ \begin{array}{l} iu_t + \Delta u   + N(u)
=i Gh = ia(x)h(x,t), \ x\in \TN, \ t>0,
\\ \\ u(x,0) =\phi (x), \end{array} \right.
\label{4.1}
\end{equation}
where $a\in C^\infty (\T ^n)$, and the nonlinearity $N(u)$ reads
\begin{equation}
\label{Nu}
N(u)=\lambda u^{\alpha _1} \overline{u}^{\alpha _2},\qquad
\alpha _1 + \alpha _2 =: \alpha + 1\ge 2, 
\end{equation}
with $\lambda \in \R$, and $\alpha , \alpha _1, \alpha _2\in \N$.
Note that for any $\alpha =2\beta\in 2\N ^*$, 
$|u|^{\alpha} u = u^{\beta +1} \overline{u}^{\beta}$. 

We introduce the number
\begin{equation}
\label{salphan}
s_{\alpha ,n}=
\left\{
\begin{array}{ll}
\frac{n}{2}- 1 \quad                            &\text{if} \ \alpha =1,\\
\frac{n}{2} -\frac{3}{4} -\frac{1}{4(n-1)}\quad &\text{if} \ \alpha =2,\\
\frac{n}{2} -\frac{2}{\alpha}                   &\text{if} \ \alpha \ge 3.
\end{array}
\right.
\end{equation}

Thus $s_{\alpha,n}=s_c:=\frac{n}{2}-\frac{2}{\alpha}$ 
(the critical Sobolev exponent obtained by scaling in NLS) for $\alpha \ge 3$, 
while $s_{\alpha, n}>s_c$ for $\alpha = 1,2$ (except for $n=\alpha =2$ where
$s_{2,2}=s_c=0$). 

By Corollary \ref{cor100} (see below), the system (\ref{4.1}) is locally
well-posed in the space $H^s (\TN)$ for $\alpha \ge 1$ and 
$s>s_{\alpha , n}$ with
$\phi \in H^s (\TN)$ and $h\in L^2_{loc} (\R , H^s (\TN))$. 

Our main concern is its exact controllability in the space 
$H^s (\TN)$.
\begin{thm} 
\label{thm10}
For given $n\ge 2$, $\alpha _1, \alpha _2\in \N$ with
$\alpha _1+\alpha _2 =:\alpha +1\ge 2$,  and $a\not\equiv 0$,  the system
(\ref{4.1}) is locally exactly controllable in the space $H^s
(\TN )$ for any $s> s_{\alpha , n}$.  More precisely, for any given
$T>0$,  there exists a number $\delta >0$ depending on $\alpha, \ n, \  T $
and $\lambda $ such that if $\phi , \ \psi \in H^s (\TN)$
satisfy
\[ \| \phi \|_s\leq \delta , 
\qquad \| \psi  \|_s \leq
\delta, \] 
then one can choose a control input $h\in L^2 (0,T; H^s
(\TN))$ such that the system (\ref{4.1}) admits a solution 
$u\in C([0,T]; H^s (\TN ))$ satisfying
\[ u(x,0) = \phi (x), \qquad u(x,T) = \psi(x) .\]
\end{thm}

The system (\ref{4.1}) can be rewritten in its equivalent integral form
\begin{equation}
u(t) = W (t) \phi + i \int ^t_0 W (t-\tau )(N(u)(\tau )) 
d\tau + \int ^t_0 W(t-\tau )[G h]
(\tau ) d\tau .\label{y-2}
\end{equation}

To prove Theorem \ref{thm10}, a smoothing property is needed for the
operator from $f$ to $u$, where
\[ u(t) = \int ^t_0 W(t-\tau ) f(\tau) d\tau .\] 
This needed smoothing property was provided in Bourgain's work 
\cite{bourgain-1,bourgain-2} where he dealt with the Cauchy problem 
for the periodic Schr\"odinger equation.

For given $s,b\in \R$, the Bourgain space $X_{s,b}$ is the space of functions
$u:\T ^n\times \R\to \C$ for which the norm
$$||u||_{X_{s,b}}=||W(-t)u(.,t)||_{H^b_t(H^s_x)}$$
is finite. Decomposing $u$ as
$$
u(x,t)=\sum_{k \in \Z ^n}
\int_\R   {\hat u}(k,\tau )e^{i(k\cdot x + \tau t)} d\tau 
$$
we have that
$$
||u||^2_{X_{s,b}}=\sum_{k\in \Z ^n}\int _\R 
\la \tau + |k|^2 \ra ^{2b}
\la k\ra ^{2s} |\hat{u}(k , \tau )|^2 d\tau
$$
where $\la y\ra :=(1+|y|^2)^{\frac12}$. 
For given $T>0$, $X_{s,b}^T$ is the restriction norm space
$$
X_{s,b}^T=\{ u_{|\T ^n\times (0,T)};\ u\in X_{s,b} \}
$$
with the restriction norm
$$
||u||_{X_{s,b}^T} = \inf \{ ||\tilde u||_{X_{s,b}}; 
\tilde u \in X_{s,b}, \ \tilde u_{\vert \TN\times (0,T)}
=u\}. 
$$

Before we proceed to show the exact controllability results, we
present the two following technical lemmas  (see e.g. \cite{tao})
which play important roles in the proof of Theorem \ref{thm10}.
\begin{lem} 
\label{lem10}
For given $T>0$ and $s,b\in \R$, there exists a 
constant $C>0$ such that
\[ \| W(t) \phi \|_{X_{s,b}^T }\leq C \| \phi \| _s \]
for any $\phi \in H^s(\T ^n).$
\end{lem}

\begin{lem} 
\label{lem11}
For given $T>0$, $b>1/2$,  and $s\in \R$, there exists a 
constant $C>0$ such that
\[ \left \| \int ^t_0 W(t-\tau) f(\tau) d \tau \right \|
_{X_{s,b}^T } \leq C \| f\| _{X_{s,b-1}^T} \] 
for any $f\in X_{s,b-1}^T$.
\end{lem}

The following multilinear estimate is crucial when applying the
contraction mapping theorem.
\begin{prop} 
\label{prop12}
Let $n\geq 2$, $\alpha \in \N ^*$ and $s > s_{\alpha,n}$. 
Then there exist some numbers $b\in (0,\frac{1}{2})$ and  $C>0$ such that
\be
\label{multilinear}
\| \prod_{i=1}^{\alpha +1} \tilde{u}_i \| _{X_{s,-b}}\leq 
C \prod_{i=1}^{\alpha +1} \| u_i\| _{X_{s,b}}\quad 
\forall u_1 ,..., u_{\alpha +1} \in X_{s,b},
\ee
where $\tilde{u}_i$ denotes $u_i$ or $\overline{u_i}$. 
\end{prop}

\begin{cor}
 \label{cor100}
Let $n\ge 2$, $\alpha \in \N ^*$, and $s>s_{\alpha ,n}$. Pick 
$u_0\in H^s(\T ^n)$ and $h\in X_{s,0}=L^2(\R; H^s(\T ^n))$. Then 
there exist two numbers $b>\frac{1}{2}$ and 
$T=T(||u_0||_{H^s(\T ^n)},||h||_{X_{s,0}})$ so that the initial-value problem
\eqref{4.1} admits a unique solution $u\in X_{s,b}^T$. 
\end{cor}
\begin{rem}
Proposition \ref{prop12}, which is proved in Appendix for the sake
of completeness, is essentially due to Bourgain. 
It was proved in
\cite{bourgain-2} when $\alpha =n=2$, and in \cite{bourgain-1} 
in Besov-type spaces when $s>s_b$, where 
\begin{equation}
\label{sbourgain}
s_b=
\left\{
\begin{array}{ll}
 s_c \quad &\text{if} \ n=2,\\
\max (s_c, \frac{3}{4}) \quad &\text{if} \ n=3,\\
\max (s_c,\frac{3n}{n+4}) \quad &\text{if}\  n\ge 4.
\end{array}
\right.
\end{equation}
Notice that $s_b>s_c$ only for $(\alpha ,n )\in \{ (2,3),(2,4),(2,5),(3,4)\}$.
The corresponding values of $s_b,s_c$ and $s_{\alpha, n}$ are reported in
Table \ref{table1}. 
\begin{table}[htbp]
\begin{center}
\begin{tabular}{|l|c|c|c|c|}
\hline
&&&&\\
$(\alpha ,n)$ & $(2,3)$       & $(2,4)$       & $(2,5)$         & $(3,4)$ \\
&&&&\\
\hline
&&&&\\
$s_b$    & $\frac{3}{4}$ & $\frac{3}{2}$ & $\frac{5}{3}$  & $\frac{3}{2}$ \\
&&&&\\
\hline
&&&&\\
$s_{\alpha,n}$ & $\frac{5}{8}$& $\frac{7}{6}$& $\frac{27}{16}$ & 
$\frac{4}{3}$\\
&&&&\\
\hline
&&&&\\
$s_c$    & $\frac{1}{2}$ & 1             & $\frac{3}{2}$   & $\frac{4}{3}$ \\
&&&&\\
\hline
\end{tabular}
\end{center}
\caption{$s_b$, $s_{\alpha, n}$ and $s_c$ for 
$(\alpha, n)\in \{(2,3),(2,4),(2,5),(3,4)\}$}
\label{table1}
\end{table}
On the other hand, $s_b=s_c < s_{\alpha, n}$ for $\alpha =2$ and $n\ge 6$. 
Sharp results for the local well-posedness of NLS on $\T ^n$ are also given
in \cite{KPV96} for $\alpha =n=1$, and
in \cite{grunrock01} for $(\alpha _1,\alpha _2)=(0,2)$ and
$2\le n\le 4$.
\end{rem}

It follows at once from Proposition \ref{prop12} that for any $T>0$, 
$s>s_{\alpha, n}$, and some
$b>1/2$, $b'>b-1$ we have 
\[
\| N(v) - N(w) \|_{X^T_{s,b'}} 
\le C (||v||^\alpha _{X^T_{s,b}} +||w||^\alpha _{X^T_{s,b}})
||v-w||_{X^T_{s,b}}\quad 
\forall v,w\in X^T_{s,b}.
\]
 
We are now in a position to give a proof of Theorem \ref{thm10}.

\medskip
 \noindent
 {\bf Proof of Theorem \ref{thm10}:} Set
\[ \omega (v, T)= i\int _0^T W (T-\tau ) N(v) (\tau)
d\tau .\]
By Theorem \ref{thm1}, if  we choose
\[ h= \Phi (\phi, \psi -\omega (v, T)),\]
then
\begin{eqnarray*}
& & \qquad  \quad W (t) \phi + \int ^t_0 W
(t-\tau )\left (iN(v)+ G \Phi (\phi, \psi -\omega (v,T)\right )
(\tau ) d\tau\\ \\
&  = &  
\left\{
\begin{array}{ll} 
\phi (x) \ \hbox{in} \ \T ^n \ & \hbox{when} \ t= 0; \\ \\
\psi (x) -\omega (v, T)+\omega (v,T)=\psi (x) \ \hbox{in} \ \T ^n , 
\ & \hbox{when} \ t=T
.\end{array} \right. \end{eqnarray*}

It suggests us to consider the nonlinear map:
\[ \Gamma (v) = W (t) \phi + i \int ^t_0 W
(t-\tau )\left (iN(v)+ G \Phi (\phi, \psi -\omega (v,T)\right )(\tau )
d\tau . \] 
The proof would be complete if we can show that this map
$\Gamma $ has a fixed point in the space $X_{s,b}^T$, with 
$b\in (\frac{1}{2},1)$. 

To this end, note that by using Lemma \ref{lem10}, Lemma \ref{lem11} and
Proposition \ref{prop12},  there
exist a number $b\in (\frac{1}{2},1)$ and some  
constants $C_j$, $j=1,2,3$ such that
\[ 
\| \Gamma (v)\| _{X_{s,b}^T} \leq C_ 1 \left (\| \phi \| _s 
+\| \psi \| _s
+\| \omega (v, T)\| _s \right ) 
+C_2 \| v\|_{X_{s,b}^T}^{\alpha +1}
\]
for any $v\in X_{s,b}^T$ and
\[
\| \Gamma (v_1)-\Gamma (v_2)\| _{X_{s,b}^T} 
\leq C_1 \| \omega (v_1, T)-\omega (v_2, T)\| _s 
+  C_3 \left ( \|v_1\|
_{X_{s,b}^T}^{\alpha} + \|v_2\| _{X_{s,b}^T}^{\alpha}\right ) \|
v_1-v_2\|_{X_{s,b}^T} \] for any $v_1,v_2 \in X_{s,b}^T$. 
Note that there exists a constant $C_4>0$ such that
\begin{eqnarray*}
 \|\omega (v,T)\|_s 
&\leq &\| \int ^t_0 W(t-\tau )
N(v) (\tau ) d\tau \|_{C([0,T]; H^s (\TN ))} \\ \\
& \leq & const \|  \int ^t_0 W(t-\tau )N(v)(\tau ) d\tau
\|_{X_{s,b}^T} \\ \\
&\leq & C_4 \|v\|^{\alpha +1} _{X_{s,b}^T} .
\end{eqnarray*}
Similarly
\[ \|\omega (v_1,T)-\omega (v_2,T)\|_s \leq C_5 
\left ( \|v_1\| _{X_{s,b}^T}^{\alpha } 
+ \|v_2\| _{X_{s,b}^T}^{\alpha}\right ) \|v_1-v_2\|_{X_{s,b}^T} .\]
As a result, by increasing the constants $C_2$ and $C_3$, we obtain 
\[ \| \Gamma (v)\| _{X_{s,b}^T} 
\leq C_ 1 (\| \phi \|_s   + \| \psi  \|_s ) 
+C_2 \| v\|_{X_{s,b}^T}^{\alpha +1} \] 
for any $v\in X_{s,b}^T $ and
\[ \| \Gamma (v_1)-\Gamma (v_2)\| _{X_{s.b}^T} 
\leq   C_3 \left ( \|v_1\| _{X_{s,b}^T}^{\alpha } 
+ \|v_2\| _{X_{s,b}^T}^{\alpha} \right )
\| v_1-v_2\|_{X_{s,b}^T} \] 
for any  $v_1 , v_2 \in X_{s,b}^T$. 
Pick $\delta >0$, $\phi , \psi  \in H^s (\TN )$
with $\|\phi \|_s + \|\psi \|_s \le \delta $,  and set 
$M=2C_1 \delta$. If $\|v\|_{X_{s,b}^T}\leq M$ and
\[ \mbox{ $\|v_j\|_{X_{s,b}^T} \leq M, \ j=1,2$,}\]
then
\begin{eqnarray*}
\| \Gamma (v)\|_{X_{s,b}^T} 
&\le & C_1 \delta + C_2  M^{\alpha + 1}\\
&\le & 2C_1 \delta = M 
\end{eqnarray*} 
as long as
\[ C_2 M^{\alpha } \le  \frac{1}{2}\cdot \] 
Choose $\delta >0$ so that $M=2C_1\delta$ fulfills
\[ C_2 M^{\alpha } \le \frac12 \
\mbox{ and }\    C_3M^\alpha  \le \frac{1}{4}, \]
and let $B_{M }$ be the ball in the space $X_{s,b}^T$ 
centered at the origin of radius $M$. For given 
$\phi,\psi  \in H^s (\TN )$ with 
$ \| \phi \|_s + \| \psi \|_s \leq \delta$, we have
\[  \| \Gamma (v)\|_{X_{s,b}^T} \leq M \] 
for any $v\in B_M$ and
\[
\| \Gamma (v_1)-\Gamma (v_2)\|_{X_{s,b}^T} 
\leq  \frac12 \| v_1 -v_2\| _{X_{s,b}^T} \]  
for any $v_1,  v_2 \in B_M$. That is to say, $\Gamma
$ is a contraction in the ball $B_M$. The proof is complete.\qed

Let us now  consider the Schr\"odinger equation posed on a cube
$\Omega =(0,\pi )^n$
\begin{equation}
iu_t + \Delta u + N(u) = ia(x)h(x,t), 
\qquad x\in \Omega, \ t\in (0,T) \label{7.4}
\end{equation}
with either the homogeneous Dirichlet boundary conditions 
\begin{equation} 
u(x,t) =0 
\qquad (x,t) \in \partial \Omega \times (0,T) \label{7.5}
\end{equation}
or the homogeneous Neumann boundary conditions
\begin{equation} 
\frac{\partial u}{\partial \nu} (x,t) = 0
\qquad (x,t) \in \partial \Omega \times (0,T). \label{7.6}
\end{equation}
The nonlinearity $N(u)$ is still as in \eqref{Nu}.

It is remarkable that internal control results with Dirichlet
(resp. Neumann) homogeneous boundary conditions can be deduced from those
already proved for periodic boundary conditions.

\begin{cor} 
\label{dirichlethom}
For given $n\ge 2$, $\alpha _1, \alpha _2\in \N$ with
$\alpha _1+\alpha _2 =:\alpha +1\ge 2$ and $\alpha$ even,  
and $a\not\equiv 0$,  the system
(\ref{7.4})-(\ref{7.5})   is locally exactly controllable in the space $H^s_D
(\Omega )$ for any $s> s_{\alpha , n}$.  More precisely, for any given
$T>0$,  there exists a number $\delta >0$ depending on $\alpha, \ n, \  T $
and $\lambda $ such that if $\phi , \ \psi \in H^s_D (\Omega )$
satisfy
\[ \| \phi \|_{H^s_D(\Omega )}\leq \delta , 
\qquad \| \psi  \|_{H^s_D(\Omega )} \leq
\delta, \] 
then one can choose a control input 
$h\in L^2 (0,T; H^s_D (\Omega ))$ such that the system 
(\ref{7.4})-(\ref{7.5}) admits a solution 
$u\in C([0,T]; H^s_D (\Omega ))$ satisfying
\[ u(x,0) = \phi (x), \qquad u(x,T) = \psi(x) .\]
\end{cor}
\begin{cor} 
\label{neumannhom}
For given $n\ge 2$, $\alpha _1, \alpha _2\in \N$ with
$\alpha _1+\alpha _2 =:\alpha +1\ge 2$  
and $a\not\equiv 0$,  the system
(\ref{7.4})-(\ref{7.6})  is locally exactly controllable in the space $H^s_N
(\Omega )$ for any $s> s_{\alpha , n}$.  More precisely, for any given
$T>0$,  there exists a number $\delta >0$ depending on $\alpha, \ n, \  T $
and $\lambda $ such that if $\phi , \ \psi \in H^s_N (\Omega )$
satisfy
\[ \| \phi \|_{H^s_N(\Omega )}\leq \delta , 
\qquad \| \psi  \|_{H^s_N(\Omega )} \leq
\delta, \] 
then one can choose a control input 
$h\in L^2 (0,T; H^s_N (\Omega ))$ such that the system 
(\ref{7.4})-(\ref{7.6}) admits a solution $u\in C([0,T]; H^s_N (\Omega ))$ 
satisfying
\[ u(x,0) = \phi (x), \qquad u(x,T) = \psi(x) .\]
\end{cor}
We shall say that a function from $(-\pi,\pi)^n$ to $\C$ is {\em odd}
(resp. {\em even}), if it is odd with respect to each coordinate $x_i$, 
$1\le i\le n$. 
The proof relies on the basic, but crucial observation that the functions in
$H^s_D(\Omega )$ (resp. $H^s_N(\Omega )$) coincide with the restrictions
to $\Omega$ of the functions in $H^s(\T ^n)$ which are odd (resp. even). 
The issue is therefore reduced to an extension of Theorem \ref{thm10}
in the framework of odd (resp. even) functions in $H^s(\T ^n)$.
Extending the function $a$ in \eqref{7.4} 
to $\T ^n$ as an even function, we notice that
the control input $h$ in Theorem \ref{thm1} can be chosen odd (resp. even)
if the functions $\phi,\psi$ are odd (resp. even). Indeed, the observability
inequality holds as well in the subspaces 
\begin{eqnarray*}
H^s_{odd}(\T ^n)  &=& \{ u\in H^s_p (\T ^n);\ 
u(x_1,...,x_{i-1},-x_i,x_{i+1},...,x_n)=-u(x)\quad\forall x\in \T ^n,
\ \forall i\},\\
H^s_{even}(\T ^n) &=& \{ u\in H^s_p (\T ^n);\ 
u(x_1,...,x_{i-1},-x_i,x_{i+1},...,x_n)=u(x)\quad\forall x\in \T ^n,
\ \forall i\}
\end{eqnarray*}
of $H^s(\T ^n)$ for $s\le 0$. On the other hand, since $u$ and $N(u)$ are
simultaneously odd (resp. even), we see that the 
contraction mapping theorem can
be applied in a space of odd (resp. even) trajectories to derive
the result in Corollary \ref{dirichlethom}
(resp. \ref{neumannhom}). Full details are provided in 
\cite{RZ2007b} for $n=1$. 

\subsection{Boundary control}

In this section we consider the Schr\"odinger equation posed on a rectangle
$\Omega =(0,l_1) \times \cdots \times (0,l_n)$
\begin{equation}
iu_t + \Delta u + N(u) = 0, 
\qquad x\in \Omega, \ t\in (0,T) \label{4.4}
\end{equation}
with either the Dirichlet boundary conditions 
\begin{equation} 
u(x,t) =1_{\Gamma _0} h(x,t) 
\qquad (x,t) \in \partial \Omega \times (0,T) \label{4.5}
\end{equation}
or the Neumann boundary conditions
\begin{equation} 
\frac{\partial u}{\partial \nu} (x,t) = 1_{\Gamma _0} h(x,t) 
\qquad (x,t) \in \partial \Omega \times (0,T). \label{4.6}
\end{equation}
When we shall consider a smooth Dirichlet controller $g$, then the boundary condition \eqref{4.5} will be replaced by 
\begin{equation}
\label{69bis}
u(x,t) = g(x) h(x,t) 
\qquad (x,t) \in \partial \Omega \times (0,T).
\end{equation}

$N(u)$ still stands for the nonlinear term in NLS. We first give a result 
(with a small control region) providing precise 
informations on the smoothness of the control input
and of the trajectories when $N(u)$ is {\em weakly} nonlinear.
%
To simplify the exposition, we assume here that
$$\Omega=(0,\pi )^n .$$
We denote by $u=W_D(t)u_0$ the solution of \eqref{B0} for $h=0$.
For given $s,b\in \R$, $X_{s,b}(\Omega )$ denotes the Bourgain 
space of functions
$u:\Omega \times \R\to \C$ for which the norm
$$||u||_{X_{s,b}(\Omega ) }= c ||W_D(-t) u(.,t)||_{H^b(\R ; 
H^s_D (\Omega ))}$$
is finite. Decomposing $u$ as
$$
u(x,t)=\sum_{k \in (\N ^*)^n}
\int_\R  {\hat u}(k,\tau )e^{i\tau t}\sin (k_1x_1)\cdots \sin (k_nx_n) d\tau 
$$
we can choose the constant $c$ so that
$$
||u||^2_{X_{s,b}(\Omega )}=\sum_{k\in (\N ^*) ^n}\int _\R 
\la \tau + |k|^2 \ra ^{2b}
\la k\ra ^{2s} |\hat{u}(k , \tau )|^2 d\tau <\infty .
$$
The restriction norm space $X_{s,b}^T(\Omega )$ is defined in the usual way
(see above the definition of $X^T_{s,b}$).
For $u\in H^s_D(\Omega )$ given, we denote by $\tilde u$ its odd extension
to $\TN =(-\pi , \pi )^n$; i.e., ${\tilde u}_{|(0,\pi )^n}=u$, and 
$\tilde u$ is odd with 
respect to each coordinate $x_i$. Note that  $\tilde u\in H^s(\T ^n)$
and $||\tilde u||_s\sim ||u||_{H^s_D(\Omega )}$. 
Defining $\tilde u(.,t)$ from $u(.,t)$ in a similar way, we observe that 
$$
||\tilde u||_{X_{s,b}^T}\sim ||u||_{X^T_{s,b}(\Omega )}.
$$
It is then clear that Lemmas \ref{lem10} and \ref{lem11} hold true with 
$W_D(t)$, $H^s_D(\Omega )$ and $X^T_{s,b}(\Omega )$ substituted to 
$W(t)$, $H^s(\TN )$ and $X^T_{s,b}$, respectively. 
We shall assume that the nonlinear term $N(u)$ satisfies the following
multilinear estimate
\be
\label{multilinearnu}
||N(u)-N(v)||_{X_{s,b'}(\Omega )} \le
c(u,v)\, ||u-v||_{X_{s,b}(\Omega )} 
\ee
where $s\in \R$, $-1/2 <b'<b\le b'+1$ and  $c(u,v)\to 0$ as $u\to 0,\ v\to 0$ in $X_{s,b}(\Omega )$.

Theorem \ref{thm3} can be extended to a semilinear context as follows.
\begin{thm}
\label{thm11}
Let $g$ be a smooth Dirichlet controller,
and let the nonlinearity $N(u)$ satisfy \eqref{Nu} and 
\eqref{multilinearnu} with 
$s\in [-1,\frac{1}{2})$, $b>0$ and $s+2b<\frac{1}{2}$. 
Pick any $T>0$. Then there exists $\delta >0$ such that for any 
$u_0,u_T \in H^s_D(\Omega )$ satisfying
$$
||u_0 ||_{H_D^s(\Omega )}\le \delta,\quad 
||u_T ||_{H_D^s(\Omega )}\le \delta
$$
one may find a control input
$h\in H^{\frac{s+1}{2}}(\T ; L^2 (\partial \Omega ))$ and a solution
$u\in C([0,T];H^s_D(\Omega )) \cap X^T_{s,b}$ of \eqref{4.4} and \eqref{69bis} such that $u(0)=u_0$ and $u(T)=u_T$.
\end{thm}
\noindent
{\em Proof. }  
For $u_T\in H^{s}_D(\Omega )$, let $h$ be the control
given by HUM which steers \eqref{AA2} from $0$ to $u_T$, namely
$h=\partial v/\partial \nu$ with $v=W_D(t)v_0$  and $v_0=S^{-1}u_T\in H^{s+2}_D(\Omega )$ (cf. Theorem \ref{thm3}). Recall that 
$h\in H^{\frac{s+1}{2}}(\T; L^2(\partial \Omega ))$ by \eqref{h1992}.
We set $u=\Lambda u_T=\Gamma S^{-1}u_T$. 
The regularity of $u$ is depicted in the following proposition.
\begin{prop}
\label{prop15} Assume that $-1\le s<1/2$ and $s+2b<1/2$. Then $\Lambda$ maps continuously $H^s_D(\Omega )$ into $C([0,T];H^s_D(\Omega ))\cap 
X_{s,b}^T(\Omega )$. 
\end{prop}
\noindent
{\em Proof of Proposition \ref{prop15}.} It follows from Proposition \ref{regularite} and Theorem \ref{thm3} that $\Lambda$ maps  
continuously $H^s_D(\Omega )$ into $C([0,T];H^s_D(\Omega ))$. Let us turn our attention to the Bourgain space $X_{s,b}^T(\Omega )$. \\
{\em Step 1.} We prove several claims used thereafter.\\
{\sc Claim 3.} For any $\gamma >1/2$, it holds
$$
\sup_{\lambda \in \R} \sum_{k\in \Z} \langle \lambda ^2 -k^2\rangle ^{-\gamma} 
<\infty.
$$ 
In what follows, $C$ denotes a constant independent of $\lambda$ and $k$ which
may vary from line to line. Pick $\lambda \in \R ^+$. For $0\le \lambda \le 1$
$$
\langle \lambda ^2 - k^2 \rangle ^{-\gamma} 
\le \langle k^2\rangle ^{-\gamma} + \langle 1-k^2\rangle ^{-\gamma} 
$$
and the result is then obvious. For $\lambda > 1$, we have 
\begin{eqnarray*}
\sum_{k\in \Z}\langle \lambda ^2 -k^2 \rangle ^{-\gamma}
&\le & C\left( 
\int_0^{\lambda -1} |\lambda ^2 -x^2|^{-\gamma }dx + 
\int_{\lambda +1}^\infty |x^2-\lambda ^2|^{-\gamma} dx +1 \right) \\
&=& C \lambda ^{1-2\gamma} \left( \int_0^{1-\lambda ^{-1}} |1-y^2|^{-\gamma}
dy + \int_{1+\lambda ^{-1}}^{+\infty} |y^2-1|^{-\gamma} dy +1\right)  \\
&\le& C \lambda ^{1-2\gamma} \left( \int_0^{1-\lambda ^{-1}} |1-y|^{-\gamma}
dy + \int_{1+\lambda ^{-1}}^{2} |y-1|^{-\gamma} dy +1\right)  \\
&\le & 
\left\{ 
\begin{array}{ll}
C\lambda ^{1-2\gamma}(\lambda ^{-1+\gamma} +1) &\text{ if } \gamma \ne 1;\\
C\lambda ^{-1}(\ln \lambda + 1) &\text{ if } \gamma =1
\end{array}
\right.
\end{eqnarray*}
and the claim follows. \\
{\sc Claim 4.} If $s\ge -1$, $0<\delta <1$, $s+2\delta <1/2$, and 
$k>1+2(s+1)$, then for some constant $C>0$
$$
S(p) := \sum_{q;|q|\ne |p|} 
\frac{q_n^{2s+2}}{||q|^2-|p|^2|^{2(1-\delta )}}
\prod_{j=1}^{n-1}\langle p_j-q_j\rangle ^{-k} \le C \langle p\rangle ^{2s+2}.
$$
Write $S(p) = S^1(p) + S^2(p)$, where the sum $S^1(p)$ is restricted to 
the $q=(q',q_n)$ with $|q'|\ge |p|$ and $|q|\ne |p|$. Noticing that 
$|q|^2-|p|^2 = q_n^2 +|q'|^2 - |p|^2 \ge q_n^2$ for such $q$, we obtain that 
$$
S^1(p) \le \sum_{q_n} q_n^{2s+4\delta -2} \sum_{q'}\prod _{j=1}^{n-1}
\langle p_j-q_j \rangle ^{-k} \le C\le C\langle p\rangle ^{2s+2} 
$$
To bound $S^2(p)$, we fix any $q'\in (\N ^*)^{n-1}$ with 
$|q'|<|p|$ and set 
$$
\lambda = \sqrt{|p|^2-|q'|^2} \ge 1.
$$
We have that 
\begin{eqnarray*}
\sum_{q_n;|q_n^2-\lambda ^2|\ge 1}
\frac{q_n^{2s+2}}{|q_n^2-\lambda ^2|^{2(1-\delta )}}
&\le& C\left( 
\int_{|x^2-\lambda ^2|\ge 1} 
\frac{x^{2s+2}}{|x^2-\lambda ^2|^{2(1-\delta )}} dx + \lambda ^{2s+2}
\right) \\
&\le& C\left( \lambda ^{2s+4\delta -1} 
\int_{|y^2-1|\ge \lambda ^{-2}}
\frac{y^{2s+2}}{|y^2-1|^{2(1-\delta)}}dy + \lambda ^{2s+2}
\right)  \\
&\le& C(\lambda ^{2s+4\delta -1}\cdot \lambda ^{2-4\delta}\cdot \ln \lambda + 
\lambda ^{2s+2})\\
&\le& C \big( p_n^{2s+2} + \sum_{j=1}^{n-1}\langle p_j^2 -q_j^2
\rangle ^{s+1}\big). 
\end{eqnarray*}
It follows that 
\begin{eqnarray*}
S^2(p) &\le & 
C\sum_{q'}(p_n^{2s+2} 
+\sum_{j=1}^{n-1}
\langle p_j^2-q_j^2 \rangle ^{s+1}) \prod_{l=1}^{n-1}
\langle p_l-q_l\rangle ^{-k} \\
&\le& C \big( p_n^{2s+2} + \sum_{j=1}^{n-1}\sum_{q_j\ge 1}
\langle p_j^2-q_j^2\rangle ^{s+1} \langle p_j-q_j\rangle ^{-k}  \big) \\
&\le& C\big( p_n^{2s+2} +\sum_{j=1}^{n-1}\sum_{q_j\ge 1}
\langle p_j+q_j \rangle ^{s+1} \langle q_j-p_j \rangle ^{-(k-s-1)} \big).
\end{eqnarray*}
To complete the proof of Claim 4, we need the following\\
{\sc Claim 5.} Let $\sigma \ge 0$ and $k>\sigma + 1$. Then there 
exists a constant $C>0$ such that 
$$
\sum_{m\ge 1} \langle m + n \rangle ^\sigma  \langle m - n \rangle  ^{-k} 
\le C n^\sigma \qquad \forall n\ge 1.
$$ 
Split the sum into $\Sigma _1 + \Sigma _2$ where 
$\Sigma _1 = \sum_{1\le m\le 3n} \langle m+n \rangle ^\sigma
\langle m-n\rangle ^{-k}$. Note that 
$$
\Sigma _1 \le \langle 4n \rangle ^\sigma \sum_{l\in \Z} 
\langle l\rangle ^{-k} \le C \langle n\rangle ^\sigma 
$$  
since $k>1$. On the other hand, noticing that 
$m-n>(m+n)/2$ for $m>3n$, we have that 
$$
\Sigma_2 \le \sum_{m>3n} \langle 2(m-n)\rangle ^\sigma 
\langle m-n\rangle ^{-k}
\le C\sum_{m>3n}\langle m-n\rangle ^{-(k-\sigma)} \le C. 
$$
Claim 5 is proved. 
Pick $k>1 + 2(s+1)\ge 1$. It follows from Claim 5 that 
$$
\sum_{q_j} \langle p_j+q_j\rangle ^{s+1}\langle p_j-q_j\rangle ^{-(k-s-1)}
\le C p_j^{s+1}.
$$
Since $s+1\ge 0$ and $p_j\ge 1$, we conclude that 
$$
S(p)\le C (p_n^{2s+2} + \langle p'\rangle ^{s+1} ) 
\le C \langle p\rangle ^{2s+2}.
$$ 
This completes the proof of Claim 4. \\
{\em Step 2.} Assume that $s<0$ and $s+2b<1/2$, and pick 
any $u_T\in H^s_D(\Omega )$ and any $\eta \in C_0^\infty (\R )$ with
$\eta (t)=1$ for $0\le t \le T$. Let $v_0=S^{-1}u_T\in H^{s+2}_D(\Omega )$ be decomposed as in \eqref{XY1}. Let us prove that $u=\Lambda u_T\in X_{s,b}^T$. 
It is sufficient to prove that 
$$
||\eta (t) u||_{X_{s,b}} \le C ||v_0||_{H^{s+2}_D(\Omega )}.
$$
Recall that $u$ is given by \eqref{XY5}-\eqref{XY6}, and that $u(t)$ may be defined
this way for all $t\in \R$ . Again, we can limit 
ourselves to proving that $u_{F_0}\in X_{s,b}^T$, where $u_{F_0}$ is the contribution due to $F_0=\{ x\in \partial \Omega; \  x_n=0 \}$ in $u$. 
$u_{F_0}$ is decomposed as  
$$
u_{F_0} =\sum_{q\in (\N ^*)^n} u_q(t) \sin (q_1x_1)\cdots \sin (q_nx_n)
$$
where 
$$u_q(t)=-\left( \frac{2}{\pi} \right) ^n 
\sum_{p\in (\N ^*)^n}
v_p\frac{e^{-i|p|^2t} - e^{-i|q|^2t} }{|q|^2-|p|^2}I_{F_0}$$
with the convention \eqref{XY4}. 
$\hat .$ denoting time Fourier transform, an application of the elementary
property 
$$\widehat{e^{irt}\eta (t)}(\tau )=\hat \eta (\tau -r)$$
yields
$$
\widehat{\eta u_q}(\tau ) = -\left( \frac{2}{\pi}\right) ^n
\left( \sum_{p;|p|\ne |q|}
v_p \frac{\hat\eta (\tau +|p|^2) - \hat\eta (\tau + |q|^2)}{|q|^2-|p|^2} I_{F_0}
+\sum_{p;|p|=|q|}
iv_p \widehat{t\eta (t) }(\tau + |q|^2)I_{F_0} \right) .
$$
For a function $w$ decomposed as 
$$
w(x,t)=\sum_{q\in (\N ^*)^n} w_q(t)\sin (q_1x_1)\cdots \sin (q_nx_n)
$$
we recall that 
$$
||w||^2_{X_{s,b}(\Omega )} 
=\sum_{q\in (\N ^*)^n}\int d\tau  \la \tau +|q|^2\ra ^{2b} \la q\ra ^{2s}
| \hat{w}_q (\tau )|^2 
$$
Therefore, it is sufficient to check that
$$
I:=\sum_{q\in (\N ^*)^n}\int d\tau  \la q \ra ^{2s}
\la \tau  + |q|^2 \ra ^{2b} | \widehat{\eta u_q}(\tau )|^2 
\le c\sum_p \la p\ra ^{2s+4} |v_p|^2.
$$
Using \eqref{XYZ}, we may write
$$
I\le c(I_1+I_2+I_3)
$$
where
\begin{eqnarray*}
I_1 &=& \sum_q\int d\tau \la q\ra ^{2s} \la \tau + |q|^2 \ra ^{2b}
\left(  \sum_{p;|p|=|q|} |v_p \, \widehat{t\eta (t)} (\tau + |q|^2)| 
p_nq_n \prod_{j=1}^{n-1} \la p_j-q_j\ra ^{-k} \right) ^2 \\
I_2 &=& \sum_q\int d\tau \la q\ra ^{2s} \la \tau + |q|^2\ra  ^{2b}
\left(  \sum_{p;|p|\ne |q|} \left\vert v_p 
\frac{\hat\eta (\tau + |q|^2)}{|q|^2-|p|^2} \right\vert 
p_nq_n \prod_{j=1}^{n-1} \la p_j-q_j\ra ^{-k} \right) ^2 \\
I_3 &=& \sum_q\int d\tau \la q\ra ^{2s} \la \tau + |q|^2\ra  ^{2b}
\left(  \sum_{p;|p|\ne |q|} \left\vert v_p 
\frac{\hat\eta (\tau + |p|^2)}{|q|^2-|p|^2} \right\vert 
p_nq_n \prod_{j=1}^{n-1} \la p_j-q_j\ra ^{-k} \right) ^2 
\end{eqnarray*}
We bound separately $I_1$, $I_2$ and $I_3$. \\
1. 
\begin{eqnarray*}
I_1 &\le& C (\int d\sigma \la \sigma \ra ^{2b} 
|\widehat{t\eta (t)} (\sigma)|^2)
\sum _q \la q\ra ^{2s} q_n^2 
\left( \sum_{p;|p|=|q|}
|v_p| p_n\prod_{j=1}^{n-1}\la p_j-q_j \ra  ^{-k} \right) ^2 \\ 
&\le& C\sum _q \la q\ra ^{2s} q_n^2 
\left( \sum_{p;|p|=|q|} |v_p|^2 p_n^2
\prod_{j=1}^{n-1}\la p_j-q_j \ra ^{-k} \right) 
\left( \sum_{p;|p|=|q|}\prod _{j=1}^{n-1}
\la p_j - q_j\ra ^{-k} \right) 
\end{eqnarray*}
where we used successively a change of variables in the integral term,
the fact that $\eta \in {\cal S }(\R )$ and Cauchy-Schwarz inequality. 
From
\begin{equation*}
\sum_{p;|p|=|q|} \prod_{j=1}^{n-1} \langle p_j-q_j\rangle ^{-k}
\le \sum_{p_1,...,p_{n-1}} \left(
\prod_{j=1}^{n-1}\langle p_j-q_j\rangle ^{-k}
\sum_{p_n;|p|=|q|} 1 \right)
\le \prod_{j=1}^{n-1} \sum_{p_j\in \Z} \la p_j\ra ^{-k}
<\infty
\end{equation*}
we deduce that 
\begin{eqnarray*}
I_1 &\le& C\sum_p |v_p|^2 |p|^2 \sum_{q;|q|=|p|}
\la q\ra ^{2s+2} \prod_{j=1}^{n-1}\la p_j-q_j\ra ^{-k} \\
&\le& C\sum_{p} |v_p|^2 |p|^{2s+4}. 
\end{eqnarray*}
2. 
\begin{eqnarray*}
I_2 &=& C (\int d\sigma \la \sigma \ra ^{2b}  
|\hat{\eta} (\sigma)|^2)
\sum _q \la q\ra ^{2s} q_n^2 
\left( \sum_{p;|p|\ne |q|}
\left\vert \frac{v_p}{ |q|^2 - |p|^2 } \right\vert
p_n\prod_{j=1}^{n-1}\la p_j-q_j \ra  ^{-k} \right) ^2 \\ 
&\le& c\sum _q \la q\ra ^{2s} q_n^2 
\left( \sum_{p;|p|\ne |q|} \frac{|v_p|^2 p_n^2}{||q|^2-|p|^2|^{2(1-\delta )}}
\prod_{j=1}^{n-1}\la p_j-q_j \ra ^{-k} \right) 
\left( \sum_{p;|p|\ne |q|} ||q|^2-|p|^2|^{-2\delta } 
\prod _{j=1}^{n-1} \la p_j - q_j\ra ^{-k} \right) 
\end{eqnarray*}
where we used Cauchy-Schwarz inequality, and $\delta > 1/4$ was chosen 
so that $s+2\delta <1/2$. From Claim 3, we obtain that 
\begin{equation*}
\sum_{p;|p|\ne |q|} | |q|^2-|p|^2|^{-2\delta} \prod_{j=1}^{n-1} 
\langle p_j - q_j\rangle ^{-k}
\le C\sum_{p'}\prod_{j=1}^{n-1} \langle p_j-q_j\rangle ^{-k} 
\sum_{p_n;|p|\ne |q|}
\langle |q|^2-|p|^2\rangle ^{-2\delta} <const.
\end{equation*}
Therefore, since $s<0$, we see that 
$$
I_2 \le C \sum_q q_n^{2s+2} \sum_{p;|p|\ne |q|}
\frac{|v_p|^2 p_n^2}{||q|^2-|p|^2|^{2(1-\delta )}}
\sum_{j=1}^{n-1} \la p_j - q_j\ra ^{-k}
$$
and from Claim 4 
$$
I_2 \le C\sum_{p} |v_p|^2 |p|^{2s+4}. 
$$
3. From the elementary estimate 
$$
\la \tau + |q|^2 \ra  
\le c \la \tau + |p|^2 \ra \la |q|^2 - |p|^2\ra 
$$
we infer that 
\be
\label{X10}
I_3\le  C\sum_q \int d\tau \la q\ra ^{2s} |q_n|^2
\left(
\sum_{p; |p|\ne |q|} |v_p|
\frac{|\hat\eta (\tau + |p|^2)| \la \tau + |p|^2\ra ^b}
{||q|^2-|p|^2|^{1-b}}  p_n \prod_{j=1}^{n-1} \la p_j-q_j\ra ^{-k}
\right) ^2.
\ee
For any fixed $\gamma >1$, we have that for some constant $c>0$ 
$$
\la \sigma \ra ^b |\hat \eta (\sigma )| \le c \la \sigma \ra  ^{-\gamma}
\qquad \forall \sigma \in \R . 
$$
Expanding the squared term in (\ref{X10}) results in 
\begin{eqnarray*}
I_3 
&\le& C\sum_q \la q\ra ^{2s} |q_n|^2
\sum_{p;|p|\ne |q|}\
\sum_{\tilde p; |\tilde p|\ne |q|} 
\frac{|v_p| \, |v_{\tilde p}| p_n\tilde p_n}
{||q|^2 -|p|^2|^{1-b} ||q|^2 -|\tilde p|^2|^{1-b}}\\
&&\quad \times (\prod_{j=1}^{n-1}\la p_j-q_j\ra ^{-k}\la \tilde p_j- q_j\ra ^{-k})
\int d\tau \la \tau + |p|^2 \ra  ^{-\gamma }
\la \tau + |\tilde p|^2 \ra ^{-\gamma}\\
&\le& C\sum_q \la q\ra ^{2s} |q_n|^2
\sum_{p;|p|\ne |q|}\
\sum_{\tilde p; |\tilde p|\ne |q|} 
\frac{|v_p| \, |v_{\tilde p}| p_n\tilde p_n}
{||q|^2 -|p|^2|^{1-b} ||q|^2 -|\tilde p|^2|^{1-b}}\\
&&\quad \times (\prod_{j=1}^{n-1}\la p_j-q_j\ra ^{-k}\la \tilde p_j- q_j\ra ^{-k})
\la |p|^2 -|\tilde p|^2\ra ^{-\gamma} 
\end{eqnarray*}
where we used the following estimate valid for $\gamma >1$ (see e.g. \cite[Lemma 7.34]{linares-ponce})
$$
\int d\tau \la  \tau + \tau _1  \ra  ^{-\gamma}
\la \tau + \tau _2 \ra  ^{-\gamma}
\le c \la \tau _1 -\tau _2\ra ^{-\gamma}.
$$
Thus
$$
I_3 \le 
C \sum_q \la q\ra ^{2s} q_n^2 
\sum_{p;|p|\ne |q|}
\frac{|v_p|^2p_n^2}{||q|^2 -|p|^2|^{2(1-b)}}(\prod_{j=1}^{n-1} \la p_j-q_j\ra ^{-k} )
\sum_{\tilde p;|\tilde p|\ne |q|} \prod_{j=1}^{n-1}
\la \tilde p_j - q_j \ra ^{-k} \la |p|^2 - |\tilde p|^2 \ra ^{-\gamma}.
$$
Since $\gamma >1/2$, it follows from Claim 3 that 
$$
\sum_{\tilde p}\prod_{j=1}^{n-1}
\la \tilde p_j -q_j \ra ^{-k} \la |p|^2 - |\tilde p| \ra ^{-\gamma}
\le \sum_{{\tilde p}_1,...,{\tilde p}_{n-1}}\prod_{j=1}^{n-1}
\la \tilde p_j -q_j\ra ^{-k} \sum_{\tilde p_n}
\la {\tilde p}_n^2 + |\tilde p'|^2 -|p|^2\ra ^{-\gamma} <const. 
$$
Thus 
$$
I_3 \le C \sum_q\la q\ra ^{2s} q_n^2 
\sum_{p;|p|\ne |q|} \frac{|v_p|^2 p_n^2}{||q|^2-|p|^2|^{2(1-b)}}
\prod_{j=1}^{n-1}\la p_j-q_j\ra ^{-k}.
$$
Using Claim 4 and the fact that $s\in [-1,0)$, we have that 
$$
I_3 \le C \sum_p |v_p|^2 |p|^2 \sum_{q;|q|\ne |p|}
\frac{q_n ^{2s+2}}{||q|^2-|p|^2|^{2(1-b)}}
\prod_{j=1}^{n-1} \la p_j-q_j\ra^{-k} 
\le \sum_p |v_p|^2 |p|^{2s+4}.
$$
\noindent{\em Step 3.}  Assume that $s+2b<1/2$ with 
$s\in [0,1/2)$. Let  $u_T$, $v_0$, $u$ and $\eta$ be as in Step 2.
Then 
\begin{eqnarray}
||\eta (t) \Gamma v_0||_{X_{s,b}} 
&\le& C||\eta D^{s+1}\Gamma v_0||_{X_{-1,b}} \nonumber \\
&\le& C \left( ||\eta (t) \Gamma (D^{s+1}v_0)||_{X_{-1,b}}
+||\eta (t) [\Gamma , D^{s+1}] v_0||_{X_{-1,b}}\right) .  \label{AAA1}
\end{eqnarray}
According to Step 2, the first term in the r.h.s. of
\eqref{AAA1} is less than $C||D^{s+1}v_0||_1\le C||v_0||_{s+2}$, for 
$-1+2b<1/2$. The contribution due to $F_0=\{ x\in \partial \Omega ; \ x_n=0 \}$ 
in $||\eta (t) [\Gamma, D^{s+1}]v_0||^2_{-1, b}$ is estimated by 
\begin{eqnarray*}
C_{F_0} &\le & \sum_q \int d\tau 
\la q \ra ^{-2} \la \tau + |q|^2 \ra ^{2b} 
\left\vert \sum_{p;|p|\ne |q|} 
v_p\frac{|p|^{s+1}-|q|^{s+1}}{|q|^2-|p|^2}
(\hat\eta (\tau + |p|^2) -\hat\eta (\tau + |q|^2))I_{F_0} \right\vert ^2 \\
&\le& C(I_2'+I_3')  
\end{eqnarray*}
where 
\begin{eqnarray*}
I_2' &=& \sum_q\int d\tau \la q\ra ^{-2} \la \tau + |q|^2 \ra ^{2b}
\left( \sum_{p;|p|\ne |q|} 
|v_p\hat\eta (\tau + |q|^2)| \frac{|p|^s + |q|^s}{|p| + |q|} p_nq_n 
\prod_{j=1}^{n-1}\la p_j-q_j\ra ^{-k} \right) ^2 ,\\
I_3' &=& \sum_q\int d\tau \la q\ra ^{-2} \la \tau + |q|^2 \ra ^{2b}
\left( \sum_{p;|p|\ne |q|} 
|v_p\hat\eta (\tau + |p|^2)| \frac{|p|^s + |q|^s}{|p| + |q|} p_nq_n 
\prod_{j=1}^{n-1}\la p_j-q_j\ra ^{-k} \right) ^2.
\end{eqnarray*}

We bound separately $I_2'$ and $I_3'$. \\
1. We have that 
\begin{eqnarray*}
I_2' &\le & 
C \big( \int d\sigma \la \sigma \ra ^{2b} |\hat\eta (\sigma )|^2  \big)
\sum_q \la q\ra ^{-2} |q_n|^2 \left\vert \sum_{p;|p|\ne |q|}
|v_p| \frac{|p|^s + |q|^s}{|p|+|q|} p_n \prod_{j=1}^{n-1}
\la p_j-q_j\ra ^{-k}  \right\vert  ^2 \\
&\le &  C\sum_q \left\vert \sum_{p}
|v_p| \frac{|p|^s + |q|^s}{|p|+|q|} p_n \prod_{j=1}^{n-1}
\la p_j-q_j\ra ^{-k}  \right\vert  ^2 \\
&\le & C\sum_p |p|^{3+2s} |v_p|^2 \\
&\le & C||v_0||^2_{s+\frac{3}{2}}. 
\end{eqnarray*}
where we used \eqref{orion1}-\eqref{orion2}. \\
2. Doing computations similar to those performed in Step 2, we obtain that
\begin{eqnarray*}
I_3' &\le & 
C\sum_q \la q\ra ^{-2} q_n^2 \sum_{p;|p|\ne |q|}
|v_p|^2 p_n^2 \frac{|p|^{2s} + |q|^{2s}}{(|p|+|q|)^2}
\left\vert |q|^2 -|p|^2\right\vert ^{2b} \prod_{j=1}^{n-1}
\la p_j-q_j\ra ^{-k} \\
&\le& C\sum_p |v_p|^2 |p|^2 \sum_{q; |q|\ne |p|}(|p|+|q|)^{2s+4b-2}
\prod_{j=1}^{n-1} \la p_j-q_j\ra ^{-k} \\
&\le& C||v_0||^2_{1}
\end{eqnarray*}
where we used the fact that $s+2b<1/2$. Since $s+2\ge 1$, we finally have that 
$$
C_{F_0}\le C ||v_0||^2_{H^{s+2}_D(\Omega )}.
$$
This completes the proof of Proposition \ref{prop15}.\qed
We can now complete the proof of Theorem \ref{thm11}. Let 
$s,b,u_0$ and $u_T$ be as in the statement of the theorem.
Using Proposition \ref{prop15} and proceeding as in the
proof of Theorem \ref{thm10}, one can show that the map 
\begin{equation}
\label{X2}
\Gamma (v) = W_D(t)u_0 + i\int_0^t W_D(t-\tau )N(v)(\tau )\, d\tau
+\Lambda (u_T - W_D(T)u_0 -\omega (v,T)) 
\end{equation}
has a fixed-point $\Gamma (v)=v$ in some closed ball $B_M\subset X_{s,b}^T(\Omega )$
provided that 
$||u_0||_{H^s_D(\Omega )} +  ||u_T||_{H^s_D(\Omega )} $ is small enough. 
Such a trajectory $v$ fulfills all the requirements of Theorem \ref{thm11}. 
In particular, $v\in X_{s,b}^T(\Omega )\cap C([0,T]; H^s_D(\Omega))$.
The smoothness of the last term in \eqref{X2} follows from
Proposition \ref{prop15}. In (\ref{X2}), we used the notation
$$
\omega (v,T) = i\int_0^T W_D(T-\tau )N(v)(\tau )d\tau. 
$$ 
Note that $\int_0^t W_D(t-\tau )N(v)(\tau )\, d\tau\in X_{s,b'+1}^T(\Omega )
\subset C([0,T]; H^s_D(\Omega ))$, by Lemma \ref{lem11}, 
\eqref{multilinearnu}, and the 
fact that $b'> -1/2$. In particular, $\omega (v,T)\in H^s_D(\Omega )$.
The proof of Theorem \ref{thm11} is achieved.
\qed
\begin{rem}
(a) Using ideas from \cite{bourgain-1}, it is likely that Theorem 
\ref{thm11} may be applied when $n\ge 2$, $\Gamma_0$ is a neighborhood of a vertex,  and $N(u)=\lambda |u|^{\alpha}u$ with $\alpha >0 $ small enough. \\
(b) The condition $s+2b<1/2$ in Proposition \ref{prop15} is actually sharp.  
Indeed, let us take $n=1$ and pick any $p\in \N^*$ and any $\eta \in {\cal S} (\R )$ 
with $|\hat\eta (\tau )|>1$ for $-1\le \tau \le 1$.  Set 
$v_0(x)=\sin (px)$ for $x\in \Omega =(0,\pi )$. With $\Gamma _0=\{ 0\}$, we have 
that $I_{F_0}=pq$ with 
$$
\widehat{\eta u_q}(\tau ) =
\left\{
\begin{array}{ll}
-\displaystyle\frac{2i}{\pi}\widehat{t\eta (t)}(\tau + p^2)p^2 \qquad &\text{ if  } q=p;\\[3mm]
-\displaystyle\frac{2}{\pi}\displaystyle\frac{\hat\eta (\tau + p^2) -\hat\eta (\tau + q^2)}{q^2-p^2}pq
\qquad &\text{ if } q\ne p.
\end{array}
\right.
$$
Therefore 
\begin{eqnarray*}
\frac{\pi ^2}{4} ||\eta u||^2_{X_{s,b}(\Omega )}
&=& \int d\tau \sum_{q;q\ne p}\la q\ra ^{2s}
\la \tau + q^2 \ra ^{2b} 
\left\vert \frac{\hat\eta (\tau + p^2) -\hat\eta (\tau +q^2)}
{q^2-p^2}\right\vert ^2 p^2q^2 \\
&&\qquad + (\int d\tau \la \tau + p^2\ra ^{2b} |\widehat{t\eta(t)} (\tau + p^2)|^2
) \la p\ra ^{2s} p^{4}\\
&=& \int d\tau \sum_{q;q\ne p}\la q\ra ^{2s} \la \tau + q^2\ra ^{2b}
\frac{|\hat\eta (\tau + p^2)|^2}{|q^2-p^2|^2} p^2q^2 + J(p)   
\end{eqnarray*}
where $|J(p)|\le Cp^{2s+4}\le C||v_0||^2_{s+2}$, according to the 
estimations of $I_1$, $I_2$, and the fact that 
$$
\int d\tau \la \tau + q^2 \ra ^{2b} |\hat\eta (\tau + p^2)
\hat\eta (\tau + q^2)|\ d\tau \le const <\infty. 
$$
Since for $q\ne p$ 
$$
\int d\tau \la \tau + q^2 \ra ^{2b} |\hat\eta (\tau + p^2)|^2
\ge \int_{-p^2-1}^{-p^2+1}d\tau \la \tau + q^2\ra ^{2b} 
\ge C |q^2-p^2|^{2b} 
$$
we have that for $s+2b\ge 1/2$,
$$
\int d\tau \sum_{q;q\ne p}\la q\ra ^{2s} \la \tau + q^2\ra ^{2b}
\frac{|\hat\eta (\tau + p^2)|^2}{|q^2-p^2|^2} p^2q^2 
\ge Cp^2\sum_{q;q>p}|q^2-p^2|^{2b-2}\la q\ra ^{2s}q^2 =\infty,
$$
therefore $\eta u\not\in X_{s,b}(\Omega )$.
The condition $s+2b<1/2$ seems related to the 
fact that any smooth function on $\T^n$ with nonnull 
boundary values belongs to the space 
$H^s_D(\Omega )$ for $s<1/2$ only. Better results will 
probably require to consider other Bourgain spaces than $X_{s,b}(\Omega )$.
\end{rem}
\begin{cor}
\label{cor11}
Let $n=1$, $\Omega =(0,\pi)$,  $\Gamma _0= \{ 0 \}$, and let the nonlinear term $N(u)$ satisfy
$$
|N(u)-N(v)| \le C (|u|^\alpha + |v|^\alpha )|u-v|,\qquad \forall u,v\in \R. 
$$
for some $\alpha \in [0,5/4)$. Let $p=\frac{4}{3}(\alpha +1) <3$.
Then there exists a number $\delta >0$ such that for any $u_0,u_T\in L^2(\Omega )$ satisfying 
$$
||u_0||_{L^2(\Omega )}<\delta, \quad ||u_T||_{L^2(\Omega )}<\delta
$$
one may find a function $h\in H^{\frac{1}{2}}(0,T)$ and a solution 
$u\in C([0,T];L^2(\Omega ))\cap L^p(0,T;L^p(\Omega ))$ 
 of \eqref{4.4}-\eqref{4.5} such that 
$u(0)=u_0$ and $u(T)=u_T$.
\end{cor}
For instance, $N_1(u)=\lambda |u|^\alpha u$ with 
$0\le \alpha <5/4$, and $N_2(u)$ of the form \eqref{Nu} with $\alpha =1$ are concerned.\\
{\em Proof.}  From the classical Strichartz estimate (see e.g. \cite{tao})
$$
||u||_{L^4(\R ; L^4(\T ))} \le C||u||_{X_{0,\frac{3}{8}}}
$$
we obtain at once the following estimates involving the spaces 
$X^T_{s,b}(\Omega )$
\begin{eqnarray*}
||u||_{L^4(0,T;L^4(\Omega ))}     &\le& C||u||_{X_{0,\frac{3}{8}}^T(\Omega )}\\
||u||_{X^T_{0,-\frac{3}{8}}(\Omega )}  &\le& 
C||u||_{L^{\frac{4}{3}} (0,T;L^{\frac{4}{3}}(\Omega ))}.
\end{eqnarray*}
Notice that for $v\in L^p(0,T;L^p(\Omega ))$, we have that 
$$
\int_0^t W_D (t-\tau )N(v)(\tau ) d\tau 
\in X^T_{0,\frac{5}{8}}(\Omega ) 
\subset C([0,T];L^2 (\Omega )) \cap L^p(0,T;L^p(\Omega )) .
$$
Indeed, 
\begin{eqnarray*}
||\int_0^t W_D(t-\tau ) N(v)(\tau ) d\tau ||_{X^T_{0,\frac{5}{8}}
(\Omega )}
&\le& C||N(v)||_{X^T_{0,-\frac{3}{8}}(\Omega )} \\
&\le& C||N(v) ||_{L^\frac{4}{3}(0,T;L^\frac{4}{3}(\Omega ))} \\
&\le& C||v||^{\alpha +1}_{L^p(0,T;L^p(\Omega ))}<\infty \cdot 
\end{eqnarray*}
In particular, $\omega (v,T)=i\int_0^T W_D (T-\tau )N(v)(\tau )d\tau
\in L^2(\Omega )$. On the other hand, by Proposition \ref{prop15}, $\Lambda$ 
maps continuously $L^2(\Omega )$ into $C([0,T];L^2(\Omega ))
\cap X_{0,b}^T (\Omega )$ for any $b<1/4$. Interpolating between
$$
X_{0,\frac{3}{8}} \subset L^4(\R ;L^4(\T )) \quad \text{ and } \quad
X_{0,0}=L^2(\R ; L^2(\T ))
$$
we obtain that 
$$
X_{0,b}\subset L^p(\R ; L^p(\T ))\quad \text{ for } 
b=\frac{3}{2}(\frac{1}{2}-\frac{1}{p})<\frac{1}{4}\cdot
$$
Therefore 
$$
\Lambda (L^2(\Omega )) \subset C([0,T];L^2(\Omega )) 
\cap L^p(0,T;L^p(\Omega )).
$$
It follows that the map 
\begin{equation*}
\Gamma (v) = W_D(t)u_0 + i\int_0^t W_D(t-\tau )N(v)(\tau )\, d\tau
+\Lambda (u_T - W_D(T)u_0 -\omega (v,T)) 
\end{equation*}
is well defined from $L^p(0,T;L^p(\Omega ))$ into 
$C([0,T];L^2(\Omega ))\cap L^p(0,T;L^p(\Omega ))$. Using the computations above,
one readily sees that $\Gamma$ contracts in some ball 
$B_M\subset L^p(0,T;L^p(\Omega ))$, provided that 
$||u_0||_{L^2(\Omega )}+||u_T||_{L^2(\Omega )}$ is small enough.\qed
\begin{cor}
\label{cor12}
Theorem \ref{thm11} may be applied when $n=2$, $\Omega =(0, \pi )^2$, $g$
is a smooth Dirichlet controller, 
$N(u)=\overline{u}^2$, $s\in (-\frac{3}{8}, -\frac{1}{3})$, 
$b\in (\frac{3}{8}, \frac{1}{2})$ with $s+2b<\frac{1}{2}$, and $b'>-\frac{1}{2}$ is sufficiently close to 
$-\frac{1}{2}$. 
\end{cor}
Corollary \ref{cor12} is a direct consequence of Theorem \ref{thm11} and of
the following result, whose proof is postponed in Appendix.
\begin{prop}
\label{prop30}
Let $s\in (-\frac{3}{8}, -\frac{1}{3})$ and  
$b\in (\frac{3}{8}, \frac{1}{2})$. Then there exists $b'\in (-\frac{1}{2},
-\frac{5}{12})$ and $C>0$ such that
\begin{eqnarray}
||\overline{v}_1\overline{v}_2||_{X_{s,b'}(\T ^2)} 
&\le& C||v_1||_{X_{s,b}(\T ^2)} ||v_2||_{X_{s,b}(\T ^2)}, 
\qquad \forall v_1,v_2\in X_{s,b}(\T ^2), \label{PP1}\\
||\overline{u}_1\overline{u}_2||_{X_{s,b'}(\Omega )} 
&\le& C||u_1||_{X_{s,b}(\Omega )} ||u_2||_{X_{s,b}(\Omega )},
\qquad \forall u_1,u_2\in X_{s,b}(\Omega ).  \label{PP2}
\end{eqnarray}
\end{prop}


Notice that if we increase the value of $s$, the state space in which the
controllability result holds has to take
into account the fact that the value (or the normal derivative) of the
function vanishes on $\partial \Omega \setminus \Gamma _0$. 
To state a result of this kind,
we limit ourselves to the situation when $\Gamma _0$ is a side, e.g.
$$
\Gamma _0=\{ 0\} \times (0,l_2)\times \cdots \times (0,l_n).
$$
Introduce the domain 
$\tilde \Omega =(-1,l_1)\times (0,l_2)\cdots \times (0,l_n)$
and a function $a\in C_0^\infty(\tilde\Omega \setminus \overline{\Omega})$, 
and consider the internal control problem
\begin{equation}
iu_t + \Delta u + N(u) = ia(x)h(x,t) , 
\qquad x\in \tilde\Omega, \ t\in (0,T). \label{4.10}
\end{equation}
Taking the restriction to 
$\Omega \times (0,T)$ of solutions of \eqref{4.10}, we obtain 
as a corollary of Theorem \ref{thm10} that both systems 
(\ref{4.4})-(\ref{4.5}) 
and (\ref{4.4})-(\ref{4.6}) are locally exactly controllable in some
subspace of $H^s (\Omega )$ for any $s> s_{\alpha ,n}$.
\begin{cor} For given $\alpha \geq 1$, $n\geq 2$, 
$\lambda \in \R$, $s> s_{\alpha, n}$
and $T>0$, there exists a constant $\delta >0$ such that for any
$u_0 , \ u_1 \in H^s (\Omega )$ satisfying
\[ \| u_i \|_{H^s (\Omega )} \leq \delta , \ i=0,1 \] 
and 
\[
u_i=\Delta u_i=\cdots =\Delta ^p u_i=0\quad  
x\in \partial\Omega \setminus \Gamma_0,\ p\le \left[\frac{2s-1}{4}\right], 
\ i=0,1
\]
\[\text{(resp.}\qquad\qquad
\frac{\partial u_i}{\partial \nu}=
\frac{\partial \Delta u_i}{\partial \nu}=\cdots =
\frac{\partial \Delta ^p u_i}{\partial \nu}=0\quad  
x\in \partial \Omega \setminus \Gamma _0, \ p\le  \left[\frac{2s-3}{4}
\right] ,\ i=0,1),
\]
then one can choose a control input $h$ such
that system (\ref{4.4})-(\ref{4.5}) (resp. system
(\ref{4.4})-(\ref{4.6})) admits a solution 
$u\in C([0,T]; H^s (\Omega ))$ with
\[ u(x,0)= u_0 (x), \qquad u(x,T) =u_1 (x). \]
\end{cor}
\noindent
\begin{rem}
By using the same extension and restriction argument, one can derive
a local controllability result in the space $H^s(\Omega )$ when 
$s>s_{\alpha, n}$ and for any given bounded smooth set $\Omega$, provided
that the control is applied on the whole boundary 
(i.e. $\Gamma _0=\partial \Omega$). A result of this kind for which 
the critical Sobolev exponent $s=s_c=s_{2,2}=0$ is reached, is given in 
\cite{RZ2008}. 
\end{rem}

\section{Stabilization}

In this section we focus on the internal stabilization of the semilinear
Schr\"odinger equation on the torus $\TN$
\be
\label{stab}
iu_t + \Delta u + N(u)=-ia^2(x)u,\qquad x\in \TN
\ee
where $a$ is any smooth real function with $a\not\equiv 0$. 

We have the following local exponential stability result which does not 
require the Geometric Control Condition. 
\begin{thm}
Let $a\in C^\infty_0(\TN)$, $a\not \equiv 0$, and let $s>s_{\alpha ,N}$. 
Then there exist some constants $\nu$, $C$ such that every solution $u$ of 
\eqref{stab} issued from
the initial state $u_0\in H^s(\TN )$ satisfies
\be
\label{decay}
||u(t)||_s\le C e^{-\nu t} ||u_0||_s \quad \forall t\ge 0.
\ee
\end{thm}
\noindent 
{\em Proof.} We proceed as in \cite{RZ2007b}. 
The operator $A_a=i\Delta - a^2$ with domain
${\mathcal D}(A_a)=H^{s+2}(\T ^n)$ generates a continuous group
$(W_a(t))_{t\in \R}$ of operators on $H^s(\T ^n)$. The first 
step is to check that the semigroup $(W_a(t))_{t\in \R ^+}$ is 
exponentially stable in $H^s(\T ^n)$. This is done in the following
\begin{prop}
\label{prop200}
There exist positive constants $C>0$ and $\nu >0$ such that
\begin{equation}
\label{s10}
||W_a(t)u_0||_s\le Ce^{-\nu t} || u_0 ||_s\qquad \forall t\ge 0. 
\end{equation}
\end{prop}    
\noindent
{\em Proof.} When $s=0$, the exponential 
stability of $(W_a(t))_{t\in \R ^+}$ is a direct consequence of 
Theorem \ref{thm1}, according to \cite{liu}. To prove (\ref{s10})
when $s=2$, we pick any $u_0\in H^2(\TN )$ and set $v:=u_t$.
Then $v$ solves the system
\begin{equation}
\left\{
\begin{array}{l}
v_t=i\Delta v -a^2(x)v,\qquad x\in \T ^n,\\
v(x,0)=v_0(x):=i\Delta u_0  (x) -a^2(x)u_0(x).
\end{array}
\right.
\end{equation}
By the property (\ref{s10}) established when $s=0$, we have
$$
||u(t)||_0\le Ce^{-\nu t}||u_0||_0, \qquad 
||v(t)||_0\le Ce^{-\nu t}||v_0||_0. 
$$
Since $i\Delta u=v+a^2 u$, we conclude that 
$$
||u(t)||_2 \le Ce^{-\nu t}||u_0||_2\qquad \forall t\ge 0.
$$ 
An easy induction yields (\ref{s10}) for any $s\in 2\N$. 
The proposition then follows by a classical interpolation argument.  \qed 

Let us now turn our attention to the stability properties of the
nonlinear system   
\begin{eqnarray*}
u_t=A_a u +iN(u), \quad  
u(.,0)=u_0
\end{eqnarray*}
that we shall write in its integral form
\begin{equation}
\label{integral}
u(t)=W_a(t)u_0 + 
i \int_0^t W_a(t-\tau )N(u)(\tau ) d \tau. 
\end{equation}
At this point, we need to establish linear estimates 
when  $W_a$ is substituted to $W$.
\begin{lem}
\label{lem4.1}
Let $T>0$, $s\geq 0$  and $0\leq b\leq 1$ be given. Then there exists a
constant $C>0$ depending only on $T$, $s$ and $b$ such that
\[ \| W_a(t) \phi \|_{X_{s,b}^T} \leq C \| \phi \| _s \]
for any $\phi \in H^s (\T ^n) $
\end{lem}
\noindent
{\em Proof.} An application of Duhamel formula gives
\begin{equation}
\label{duhamel}
W_a(t)\phi =W(t)\phi  - \int_0^t W(t-\tau )(a^2W_a(\tau )\phi)d\tau. 
\end{equation}
It follows that
\begin{eqnarray*}
|| W_a(t)\phi ||_{X_{s,b}^T}
&\le& || W(t)\phi ||_{X_{s,b}^T} +
||\int_0^t W(t-\tau )(a^2 W_a(\tau )\phi ) d\tau||_{X_{s,b}^T} \\
&\le& C||\phi||_s + C ||a^2 W_a(t)\phi||_{ X_{s,b-1}^T} \\
&\le& C||\phi||_s + C||W_a(t)\phi||_{ L^2 (0,T;H^s(\T ^n)) }
\qquad (\hbox{as}\ b-1\le 0)\\
&\le& C||\phi||_s, 
\end{eqnarray*}
as desired. \qed 

\begin{lem}
\label{lem4.2}
Let $T>0$, $s\geq 0$, and $b\in (\frac12 , 1)$  be given.  
Then there exists a constant $C>0$ depending only on 
$T$, $s$ and $b$ such that
\[ \left \| \int ^t_0 W_a(t-\tau) f(\tau) d \tau \right \|_{X_{s,b}^T} 
\leq C\| f\| _{X_{s,b-1}^T} \] for any $f\in X_{s,b-1}^T.$
\end{lem}
\noindent
{\em Proof.}
It follows from (\ref{duhamel}) that
$$
\int_0^t W_a(t-\tau )f(\tau )d\tau  = 
\int_0^t W  (t-\tau )f(\tau )d\tau  - \int_0^t W(t-\tau ) a^2 
\left( \int_0^\tau W_a (\tau -s)f(s)ds\right) d\tau, 
$$
hence
\begin{eqnarray*}
||\int_0^t W_a(t-\tau )f(\tau )d\tau ||_{X_{s,b}^T}
&\le& C||f||_{X_{s,b-1}^T} 
+ C||a^2 \int_0^t W_a (t-s) f(s)ds||_{X_{s,b-1}^T}\\
&\le& C||f||_{X_{s,b-1}^T} 
+ C||\int_0^t W_a (t-s) f(s)ds||_{X_{s,0}^T}\\
&\le& 
C||f||_{X_{s,b-1}^T} 
+ C T^\alpha ||\int_0^t  W_a(t-s)f(s)\, ds ||_{X_{s,b}^T}
\end{eqnarray*}
for some constant $\alpha >0$, by virtue of 
Lemmas \ref{lem10} and \cite[Lemma 2.11]{tao}.
The result follows at once if $T$ is small enough, 
say $T<T_0$. For $T\ge T_0$, the result 
follows from Lemma \ref{lem4.1} and an easy induction.
\qed 

Let us now proceed to the proof of the exponential stability of the 
system (\ref{stab}). Pick a number $s\ge 0$. According to
Proposition \ref{prop200}, there exist positive constants 
$C,\nu$ such that 
$$
||W_a(t) u_0||_s \le Ce^{-\nu t } ||u_0||_s \qquad \forall t\ge 0. 
$$    
Pick a time $T>0$ such that
$$
Ce^{-\nu T}< \frac{1}{4}
$$  
and fix a number $b\in (\frac12,1)$. We seek a solution $u$ of the
integral equation (\ref{integral}) in the form of a fixed point of the map 
$$
\Gamma (u)=W_a(t)u_0 + i\int_0^t W_a(t-\tau )N(u)(\tau )d\tau
$$
in some ball $B_M$ of the space $X_{s,b}^T$. This will be done provided that
$||u_0||_s\le \delta$ where $\delta $ is a small number to be determined. 
Furthermore, to ensure the exponential stability,  
$\delta$ and $M$ will be chosen in such a way that
$||u(T)||_s\le ||u_0||_s/2$. Pick for the moment any $\delta >0$ and $M>0$, 
and let $u_0\in H^s(\T ^n)$ be such that 
$||u_0||_s \le \delta $. By computations
similar to those displayed in the proof of Theorem \ref{thm10} with $W_a(t)$
substituted to $W(t)$, we arrive to 
$$
||\Gamma (u)||_{X_{s,b}^T} \le c ||u_0||_s + cM^{\alpha +1} \qquad
\forall u\in B_M
$$ 
and 
$$
||\Gamma (u)-\Gamma (v)||_{X_{s,b}^T} \le c M^{\alpha } ||u-v||_{X_{s,b}^T}
\qquad \forall u,v \in B_M$$
for some constant $c>0$ independent of $\delta$, $M$, and $u_0$.
On the other hand, using the estimate of $||\omega (T,u)||_s$ in
the proof of Theorem \ref{thm10}, we obtain
\begin{eqnarray*}
||\Gamma (u)(T)||_s 
&\le& ||W_a(T)u_0||_s + ||\int_0^T W_a(T-t) N(u)(t)dt ||_s \\
&\le& \frac{1}{4} ||u_0||_s +cM^{\alpha +1}.
\end{eqnarray*}   
Pick $\delta =4cM^{\alpha +1}$ where $M>0$ is chosen so that 
$$
(4c^2+c)M^{\alpha +1}\le M \ \hbox{ and }  \  c M^\alpha \le \frac{1}{2}. 
$$
Then we have 
\begin{eqnarray*}
||\Gamma (u)||_{X_{s,b}^T} &\le& M \qquad \forall u\in B_M\\ 
||\Gamma (u)-\Gamma (v)||_{X_{s,b}^T} &\le& \frac{1}{2}||u-v||_{X_{s,b}^T}
\qquad \forall u,v\in B_M.
\end{eqnarray*}
Thus the map $\Gamma$, which is a contraction in $B_M$, has a
fixed point $u\in B_M$. By construction,  $u$ fulfills
$$
||u(T)||_s  = ||\Gamma (u)(T) ||_s \le \frac{\delta }{2}. 
$$
Assume now that $0<||u_0||_s<\delta$. Changing $\delta$ into 
$\delta ':=||u_0||_s$ and $M$ into 
$M':=(\delta '/\delta)^{\frac{1}{\alpha +1}}M$, we 
obtain that $||u(T)||_s\le ||u_0||_s/2$, and an obvious induction yields
$||u(kT)||_s\le 2^{-k}||u_0||_s$ for any $k\ge 0$. As 
$X_{s,b}^T\subset C([0,T];H^s(\T ^n))$ for $b>1/2$, and 
$||u||_{X_{s,b}^T}\le M=(\delta/(4c))^{\frac{1}{\alpha + 1}}$, we
infer by the semigroup property that there exist some constants
$C'>0, \nu'>0$ such that 
$$
||u(t)||_s\le C'e^{-\nu 't}||u_0||_s. 
$$  
The proof is complete. \qed

\section{Appendix}
\subsection{Proof of Proposition \ref{prop12}.}

We proceed as in \cite[pp. 115-118]{bourgain-2}. 
We first introduce some notations.
Let $|x|_\infty :=\sup_{1\le i\le n}|x_i|$ 
for $x=(x_i)_{1\le i\le n}\in \R^n$. We introduce a 
dyadic partition of $\R ^n$ 
$$
\Z ^n = \cup _{j\ge 0} D_j,
$$ 
where $D_0 =\{ 0 \}$, and $D_j=\{ k\in \Z ^n;\ 2^{j-1} \le |k|_\infty <2^j\}$
for $j\ge 1$. For any H\"older exponent $p,q\in [1,+\infty ]$, we write 
$L_t^pL_x^q$ for $L^p(\R _t, L^q (\T ^n_x))$. The (discrete) cube of center
$x_0\in \R ^n$ and sidelength $2R>0$ is 
$$
Q(x_0,R)=\{  k\in \Z ^n;\ |k-x_0|_\infty \le R \} .
$$
The Strichartz estimate (\cite{bourgain-1},\cite{grunrock00})
$$
||u||_{L^4_tL^4_x} \le c||u||_{X_{s,b}},
\qquad s>\frac{n}{2} -\frac{n+2}{4}, \  b>\frac{1}{2},
$$
when combined with the standard estimates
\begin{eqnarray*}
|| u ||_{L^\infty _t L^2_x} &\le& c ||u||_{X_{0,b}},\quad b>\frac{1}{2} \\
|| u ||_{L^2_tL^2_x} &=&||u||_{X_{0,0}}
\end{eqnarray*}
and Sobolev embedding theorem, gives by interpolation the following result.
\begin{lem}(\cite[cor. 2.2]{grunrock00}) 
\label{interpolation}
 Let $n\ge 2$.\\
 (i) For all $p,q,s$ satisfying 
\begin{equation}
\label{S20}
0< \frac{1}{p}\le \frac{1}{4},\ 0< \frac{1}{q} \le \frac{1}{2} -
\frac{1}{p},
\ s>\frac{n}{2} -\frac{2}{p} -\frac{n}{q},
\end{equation}
there exists a number $b\in (0,\frac{1}{2})$ such that for all
$u\in X_{s,b}$, it holds
\begin{equation}
\label{S21}
||u||_{L^p_tL^q_x} \le c ||u||_{X_{s,b}}
\end{equation}
(ii) For all $p,q,s,b$ satisfying
\begin{equation}
\label{S22}
0\le \frac{1}{p}\le \frac{1}{q}\le \frac{1}{2}\le
\frac{1}{p} +  \frac{1}{q} \le 1,\
s>(n-2)(\frac{1}{2} -\frac{1}{q}),\
\text{ and } b > 1 - \frac{1}{p}-\frac{1}{q}
\end{equation}
then for all $u\in X_{s,b}$, \eqref{S21} holds.
\end{lem} 
Let ${\cal F}_x$ denote the Fourier transform in $x$, and let
$1_{Q}$ denote the characteristic function of the 
cube $Q$. The following result,
inspired by an observation made in \cite{bourgain-1}, indicates that
for a function spatially supported in a cube, only the sidelength of
the cube (not its center) comes into play in \eqref{S21}. 
\begin{lem} (\cite[Lemma 2.4]{grunrock00})
Assume that for $p,q,s,b$ the estimate \eqref{S21} is valid. Then there exists
a constant $c>0$ such that for any cube
$Q$ of center $x_0\in \R ^n$ and sidelength $R>0$ it holds 
\be
||({\cal F}_x^{-1} 1_Q {\cal F}_x) u||_{L^p_t L^q_x}
\le c R^s ||u||_{X_{0,b}}\cdot
\ee
\end{lem}
It follows that if \eqref{S20} (or \eqref{S22}) holds and if $u=u(x,t)$ is
a function decomposed as 
$$
u(x,t)=\sum_{|k-x_0|_\infty \le R} \int_{\R} 
{\hat u}(k,\tau ) e^{i(k\cdot x + \tau t)} d\tau
$$
then
\be
||u||_{L^p_tL^q_x} \le c R^s ||u||_{X_{0,b}}
=cR^s \left(\sum_{|k-x_0|_\infty \le R}\int_{\R}
\la \tau + |k|^2 \ra ^{2b} |\hat u (k,\tau )|^2  d\tau
\right) ^{\frac{1}{2}}.
\label{S15}
\ee
Let the functions $u_1,...,u_{\alpha + 1}\in X_{s,b}$ be given, where 
$s$ and $b$ denote some positive numbers,  and let us set 
$$
u= \tilde{u}_1 \tilde{u}_2 \cdots \tilde{u}_{\alpha  + 1}
$$
where $\tilde{u}_i$ is $u_i$ or $\overline{u_i}$. To estimate 
$||u||_{X_{s,-b}}$ we proceed by duality, estimating the integral
$\int_{\R}\int_{\T ^n}u\overline{v}dxdt$ for any $v\in X_{-s,b}$ with
$||v||_{X_{-s,b}}\le 1$. By Plancherel theorem  
\begin{eqnarray*}
\int_{\R}\int_{\T ^n} u\overline{v}\, dxdt 
&=& \sum_{k\in \Z ^n}\int_{\R} \hat u(k,\tau) \overline{\hat v}(k,\tau )d\tau\\
&=& \sum_{k_1\cdots k_{\alpha +1}}\int_{\tau_1 \cdots \tau _{\alpha +1}}
\la k\ra ^s \big(\prod_{i=1}^{\alpha +1} \hat{\tilde u}_i(k_i,\tau_i)\big)
\la k\ra ^{-s}\overline{\hat v}(k,\tau) 
\end{eqnarray*}
where $k=k_1 + \cdots + k_{\alpha +1}$ and $\tau =\tau _1 + \cdots +
\tau _{\alpha +1}$. Notice that
$\hat{\overline u}(k_i,\tau _i)=\overline{\hat{u}_i (-k_i,-\tau _i)}$. 
Writing $k_i\in D_{j_i}$, $j_i\ge 0$, we obtain
$$
\vert \int_{\R}\int_{\T ^n} u\overline{v}\, dxdt \vert 
\le\sum_{j_1\cdots j_{\alpha +1}} \sum_{k_i\in D_{j_i}}
\int_{\tau _1\cdots \tau_{\alpha +1}}
\la k\ra ^s (\prod _{i=1}^{\alpha +1}
\vert\hat{u}_i(k_i,\tau _i)\vert ) \la k \ra ^{-s}
|\hat v(k,\tau )|,
$$
where now 
$k=\pm k_1 \cdots \pm k_{\alpha +1}$, 
$\tau =\pm \tau _1 \cdots \pm \tau_{\alpha +1}$
($+k_i$ if $\tilde {u}_i=u_i$, $-k_i$ if $\tilde{u}_i=
\overline{u_i}$, and the same for $\pm\tau _i$).
We shall
focus on the sum $\Sigma =\sum_{j_1\ge j_2\ge \cdots \ge j_{\alpha +1}}$, the
other contributions leading to similar bounds. 
As $|k_i|_\infty \le 2|k_1|_\infty$ for $i\ge 2$,
we have that 
$$
\Sigma\le c\sum_{j_1\ge \cdots\ge j_{\alpha +1}}2^{j_1s}
\sum_{k_i\in D_{j_i}}
\int_{\tau _1\cdots \tau_{\alpha +1}}
 (\prod _{i=1}^{\alpha +1}
\vert\hat{u}_i(k_i,\tau _i)\vert ) \la k \ra ^{-s}
|\hat v(k,\tau )|.
$$
Pick $\gamma \in\N ^*$ with 
$$
\alpha \le 2^{\gamma -2}
$$
and split $\Sigma$ into $\Sigma _1 + \Sigma _2$ where $\Sigma _1$ corresponds
to the $j_1,...,j_{\alpha +1}$ for which
$$
j_1\ge j_2+\gamma +2 \ge j_2\ge j_3 \ge \cdots \ge j_{\alpha +1}.
$$ 
Consider a ``partition'' of $D_{j_1}$ into a collection of cubes $Q_l$ of 
sidelength $2^{j_2}$
$$
D_{j_1}=\cup_{l} Q_l.
$$
Note that each $k\in D_{j_1}$ belongs to at most $2^n$ cubes $Q_l$.
For any $l$, we denote by ${\tilde Q}_l$ the cube of sidelength
$2^{j_2+\gamma}$ with the same center as $Q_l$ if $k=k_1\pm k_2\cdots$, and
with center the opposite of that of $Q_l$ if $k=-k_1\pm k_2\cdots$. 
We claim that $k\in {\tilde Q}_l$ when $k_1\in Q_l$ and 
$k_i\in D_{j_i}$ for $i\ge 2$. Indeed
\be
\label{S30}
|k_2|_\infty +\cdots +|k_{\alpha +1}|_\infty
\le \alpha 2^{j_2} \le 2^{j_2+\gamma -2},
\ee  
hence if $Q_l=Q(x_0,2^{j_2-1})$
$$
|\pm x_0-k|_\infty \le |\pm x_0-\pm k_1|_\infty + |k_2|_\infty +\cdots + 
|k_{\alpha +1}|_\infty 
\le 2^{j_2 -1} + 2^{j_2+\gamma -2} \le 2^{j_2+\gamma -1}.
$$
Notice also that 
${\tilde Q}_l\subset D_{j_1-1}\cup D_{j_1} \cup  D_{j_1 +1}$ since 
the sidelength of $\tilde{Q}_l$ is at most $2^{j_1-2}$ and $Q_l\subset
D_{j_1}$.  It follows that 
$$
\Sigma _1 \le c \!\!\!\!\!\!\!\!
\sum_{
\begin{array}{c}
\scriptstyle j_1\ge j_2+\gamma + 2\\
\scriptstyle j_2\ge j_3\ge \cdots \ge j_{\alpha +1}
\end{array}}
\!\!\!\!\!\!\!\!
2^{j_1s}
\sum_l
\sum_{k_1\in Q_l}
\!\!\!\!\!\!\!\!
\sum_{
\begin{array}{c}
\scriptstyle k_2\in D_{j_2},\\
\scriptstyle k_{\alpha +1}\in D_{j_{\alpha +1}}
\end{array}}
\int_{\tau _1\cdots \tau _{\alpha +1}}
(\prod_{i=1}^{\alpha +1} | \hat{u}_i(k_i,\tau _i)|)
1_{\tilde{Q}_l}(k) \la k\ra ^{-s} |\hat v(k, \tau )|.
$$
Let us introduce the functions
\begin{eqnarray*}
f_l(x,t) &=& \sum_{k\in Q_l} 
\int_\R |\hat{u}_1 (k,\tau )|e^{i(k\cdot x + \tau t)} d\tau \\
g_l(x,t) &=& \sum_{k\in \tilde{Q}_l} 
\int_\R \la k\ra ^{-s} |\hat v (k, \tau )| e^{i(k\cdot x + \tau t)} d\tau
\end{eqnarray*}
and
$$
h_i(x,t)=\sum_{k\in D_{j_i}} \int_\R |\hat{u}_i(k,\tau )|
e^{i(k\cdot x + \tau t)}d\tau \quad \text{ for } i=2,...,\alpha +1.
$$
By Plancherel theorem
$$
\Sigma_1 \le c\!\!\!\!\!\!\!
\sum_{
\begin{array}{c}
\scriptstyle j_1\ge j_2+\gamma +2\\
\scriptstyle j_2\ge j_3\ge \cdots \ge j_{\alpha +1}
\end{array}}\!\!\!\!\!\!\!
2^{j_1s}
\sum_l\int_{\R} \int_{\T ^n} 
|f_lh_2\cdots h_{\alpha +1}g_l|\, dxdt.
$$
Pick H\"older exponents $p_1,q_1,p_2,q_2\in [1,\infty )$ such that
\begin{eqnarray}
\frac{3}{p_1} + \frac{\alpha -1}{p_2} &=& 1 \label{S50} \\
\frac{3}{q_1} + \frac{\alpha -1}{q_2} &=& 1 \label{S51}
\end{eqnarray}
We have that 
$$
\int_{\R}\int_{\T ^n}|f_lh_2\cdots h_{\alpha +1}g_l| dxdt
\le 
||f_l||_{L_t^{p_1}L_x^{q_1}} 
||g_l||_{L_t^{p_1}L_x^{q_1}}
||h_2||_{L_t^{p_1}L_x^{q_1}}
\prod_{i=3}^{\alpha +1}||h_i||_{L_t^{p_2}L_x^{q_2}}.
$$
Assume that for some exponents $s_1,b_1,s_2,b_2$ the following estimates
hold
\begin{eqnarray}
||u||_{L_t^{p_1}L_x^{q_1}} &\le& c ||u||_{X_{s_1,b_1}}, \label{S40} \\
||u||_{L_t^{p_2}L_x^{q_2}} &\le& c ||u||_{X_{s_2,b_2}}. \label{S41}
\end{eqnarray}
Then, by \eqref{S15} and the fact that the sidelength of $Q_l$ (resp. 
$\tilde{Q}_l$) is $2^{j_2}$ (resp. $2^{j_2+\gamma}$), we have
\begin{eqnarray}
||f_l||_{L_t^{p_1}L_x^{q_1}} &\le& c2^{j_2s_1}
\big( \sum_{k\in Q_l}\int_\tau \la \tau + |k|^2\ra ^{2b_1}
|\hat u_1|^2 \big) ^{\frac{1}{2}} \label{S300} \\
||g_l||_{L_t^{p_1}L_x^{q_1}} &\le& c2^{j_2s_1}
\big( \sum_{k\in \tilde{Q}_l}\int_\tau \la \tau + |k|^2\ra ^{2b_1}
\la k \ra ^{-2s}
|\hat v|^2 \big) ^{\frac{1}{2}} \label{S301} \\
||h_2||_{L_t^{p_1}L_x^{q_1}} &\le& c2^{j_2s_1}
\big( \sum_{k\in D_{j_2}}\int_\tau \la \tau + |k|^2\ra ^{2b_1}
|\hat u_2|^2 \big) ^{\frac{1}{2}} \label{S302}
\end{eqnarray}
and for $i=3,...,\alpha +1$
\begin{eqnarray}
||h_i||_{L_t^{p_2}L_x^{q_2}} 
&\le& c2^{j_is_2}
\big( \sum_{k\in D_{j_i}}\int_\tau \la \tau + |k|^2\ra ^{2b_2}
|\hat u_i|^2 \big) ^{\frac{1}{2}}\nonumber \\
&\le& c\big( \sum_{k\in D_{j_i}}
\int_\tau \la \tau + |k|^2 \ra ^{2b_2} \la k \ra ^{2s_2} |\hat u_i |^2  
\big) ^{\frac{1}{2}}.
 \label{S303}
\end{eqnarray} 
Using Cauchy-Schwarz in $\sum_l$, we obtain
\begin{eqnarray*}
\Sigma _1 
&\le& c\!\!\!\!\!\!\!\!\!\!
\sum_{
\begin{array}{c}
\scriptstyle j_1\ge j_2+\gamma +2\\
\scriptstyle j_2\ge j_3\ge \cdots \ge j_{\alpha +1}
\end{array}} 
\!\!\!\!\!\!\!\!\!\!
2^{j_1s+3j_2s_1}
\big( 
\sum_l\sum_{k\in Q_l}
\int_\tau \la \tau + |k|^2\ra ^{2b_1}|\hat u_1|^2 
\big) ^{\frac{1}{2}}
\big(
\sum_l\sum_{k\in \tilde{Q}_l}
\int_\tau \la \tau + |k|^2\ra ^{2b_1}
\la k\ra ^{-2s} |\hat v|^2 \big) ^{\frac{1}{2}}\\
&& \qquad
\big( 
\sum_{k\in D_{j_2}}\int_\tau \la \tau + |k|^2\ra ^{2b_1}|\hat u_2|^2 
\big) ^{\frac{1}{2}}
\prod_{i=3}^{\alpha +1}
\big(
\sum_{k\in D_{j_i}}\int_\tau 
\la \tau + |k|^2\ra ^{2b_2}
\la k\ra ^{2s_2} |\hat u_i|^2 \big) ^{\frac{1}{2}}\\
&\le& c\!\!\!\!\!\!\!\!\!\!
\sum_{
\begin{array}{c}
\scriptstyle j_1\ge j_2+\gamma +2\\
\scriptstyle j_2\ge j_3\ge \cdots \ge j_{\alpha +1}
\end{array}} 
\!\!\!\!\!\!\!\!\!\!
\big( 
\sum_{k\in D_{j_1}}\int_\tau \la \tau + |k|^2\ra ^{2b_1}
\la k\ra ^{2s}|\hat u_1|^2 
\big) ^{\frac{1}{2}}
\big(
\sum_{k\in D_{j_1-1}\cup D_{j_1}\cup D_{j_1+1}}
\int_\tau \la \tau + |k|^2\ra ^{2b_1}
\la k\ra ^{-2s} |\hat v|^2 \big) ^{\frac{1}{2}}\\
&& \qquad
\big( 
\sum_{k\in D_{j_2}}\int_\tau \la \tau + |k|^2\ra ^{2b_1}
\la k\ra ^{6s_1} |\hat u_2|^2 \big) ^{\frac{1}{2}}
\prod_{i=3}^{\alpha +1}
\big(
\sum_{k\in D_{j_i}}\int_\tau 
\la \tau + |k|^2\ra ^{2b_2}
\la k\ra ^{2s_2} |\hat u_i|^2 \big) ^{\frac{1}{2}}.
\end{eqnarray*}
We used the fact that a point $k\in D_{j_1-1}\cup D_{j_1}\cup D_{j_1+1}$
belongs to (at most) a finite number of cubes ${\tilde Q}_l$, bounded by 
$(2^{\gamma +2}+1)^n$. 
A sum 
$\sum_{j_i\ge 0}\big( 
\sum_{k\in D_{j_i}}
\int_\tau  \la \tau + |k|^2 \ra ^{2b_2} \la k\ra ^{2s_2} 
|\hat u_i|^2\big) ^{\frac{1}{2}} $
can be estimated by $c||u_i|| _{X_{s_2+\varepsilon , b_2}}$ for any
$\varepsilon >0$ thanks to Cauchy-Schwarz. Summing successively in 
$k_{\alpha +1}, ..., k_1$, we arrive at 
$$
\Sigma _1 
\le c||u_1||_{X_{s,b_1}} ||v||_{X_{-s,b_1}}
||u_2||_{X_{3s_1+\varepsilon, b_1}}
\prod_{i=3}^{\alpha + 1}
||u_i||_{X_{s_2+\varepsilon, b_2}}.
$$ 
The same bound for $\Sigma _2$ can be obtained by a more simple analysis.
Indeed, as $j_1\le j_2+\gamma +1$ in the sum over $j_1,...,j_{\alpha +1}$,
we obtain
$$
\Sigma _2 \le 
c\!\!\!\!\!\!\!\!\!\!
\sum_{
\begin{array}{c}
\scriptstyle j_1\le  j_2+\gamma +1 \\
\scriptstyle j_2\ge j_3\ge \cdots \ge j_{\alpha +1}
\end{array}} 
\!\!\!\!\!\!\!\!\!\!
2^{j_1s}
\int_{\R} \int_{\T ^n} |f h_2\cdots h_{\alpha + 1} g| dxdt,
$$ 
where 
\begin{eqnarray*}
f(x,t) &=& \sum_{k\in D_{j_1}} 
\int_\R |\hat{u}_1 (k,\tau )|e^{i(k\cdot x + \tau t)} d\tau \\
g(x,t) &=& \sum_{|k|\le (2^{\gamma +1} + \alpha) 2^{j_2}} 
\int_\R \la k\ra ^{-s} |\hat v (k, \tau )| e^{i(k\cdot x + \tau t)} d\tau
\end{eqnarray*}
and $h_2,...,h_{\alpha +1}$ as above. Since 
$2^{j_1s_1}\le c 2^{j_2s_1}$, we still have
\begin{eqnarray*}
||f||_{L_t^{p_1}L_x^{q_1}} 
&\le& c2^{j_2s_1} \big( \sum_{k\in D_{j_1}}
\int_\tau \la \tau + |k|^2\ra ^{2b_1}
|\hat u_1|^2 \big) ^{\frac{1}{2}} \\
||g||_{L_t^{p_1}L_x^{q_1}} 
&\le& c2^{j_2s_1} \big( \sum_{k\in \Z ^n}\int_\tau \la \tau + |k|^2\ra ^{2b_1}
\la k \ra ^{-2s}
|\hat v|^2 \big) ^{\frac{1}{2}} 
\end{eqnarray*}
Next, $\Sigma _2$ is estimated as $\Sigma _1$ (see above). 
At this stage, we have
proved that 
\begin{equation}
\label{S500}
\Sigma \le 
c||u_1||_{X_{s,b_1}} ||v||_{X_{-s,b_1}}
||u_2||_{X_{3s_1+\varepsilon , b_1}} 
\prod_{i=3}^{\alpha +1} ||u_i||_{X_{s_2+\varepsilon}, b_2}
\end{equation}
where $\varepsilon >0$ is arbitrary small, the exponents $s_1,b_1,s_2,b_2$ 
are taken so that \eqref{S40}-\eqref{S41} are satisfied, with the
H\"older exponents $p_1,q_1,p_2,q_2$ satisfying \eqref{S50}-\eqref{S51}. 
The proof will be complete if, in addition, we have 
$$
s\ge \sup \{ 3s_1+\varepsilon, s_2+\varepsilon \} ,\quad
b_1<\frac{1}{2}, b_2<\frac{1}{2}.
$$
We distinguish three cases:
(i) $\alpha \ge 3$; (ii) $\alpha =2$; (iii) $\alpha =1$.\\
\underline{(i) $\alpha \ge 3$}\\
We aim to reach any value $s>s_c$. To find the sets of exponents
$(p_1,q_1,s_1,b_1)$, $(p_2,q_2,s_2,b_2)$ satisfying \eqref{S20}, 
\eqref{S50} and
 \eqref{S51}, and leading to the ``smallest'' value of $s$, we are let to
minimize the functional $\sup \{ 3\sigma _1,\sigma _2\}$, where
\begin{eqnarray}
\sigma _1 &=& \frac{n}{2} - (\frac{2}{p_1}+\frac{n}{q_1}) \label{T1}\\
\sigma _2 &=& \frac{n}{2} - (\frac{2}{p_2}+\frac{n}{q_2}) \label{T2}
\end{eqnarray}
under the constraints 
\begin{eqnarray}
&&4\le p_1<\infty \label{T3}\\
&&0 < \frac{1}{q_1}\le \frac{1}{2} - \frac{1}{p_1} \label{T30}\\
&&4\le p_2 <\infty \label{T4}\\
&&0 < \frac{1}{q_2}\le \frac{1}{2} - \frac{1}{p_2} \label{T40}\\
&&\frac{3}{p_1} +\frac{\alpha -1}{p_2} = 1 \label{T5}\\
&&\frac{3}{q_1} + \frac{\alpha -1}{q_2} = 1. \label{T6}
\end{eqnarray}
At this point, it is convenient to introduce the numbers $r_1,r_2$ with 
\begin{eqnarray}
\frac{1}{r_1} &=& \frac{2}{p_1} + \frac{n}{q_1} \label{T7}\\
\frac{1}{r_2} &=& \frac{2}{p_2} + \frac{n}{q_2}\cdot \label{T8}
\end{eqnarray}
Note that, by \eqref{T5}-\eqref{T6}, 
\be
\frac{3}{r_1} + \frac{\alpha -1}{r_2}=n+2.
\label{T9}
\ee
Therefore, 
$3\sigma _1=\frac{n}{2}-2+\frac{\alpha -1}{r_2}$ (resp. 
 $\sigma _2=\frac{n}{2}-\frac{1}{r_2}$) is a nonincreasing function 
(resp. a nondecreasing function) of $r_2$. Thus the least value of 
$\sup \{3\sigma_1, \sigma _2\}$ 
is achieved when $3\sigma _1=\sigma _2$, which yields 
\be
\label{S1000}
r_2=\frac{\alpha}{2}, \ r_1 = 3 (n+\frac{2}{\alpha})^{-1}, \quad
3\sigma _1=\sigma_2=\frac{n}{2}-\frac{2}{\alpha}\cdot
\ee
It remains to find $p_1,q_1,p_2,q_2$ satisfying \eqref{T3}-\eqref{T8}. 
Note first that \eqref{T6} is satisfied whenever \eqref{T5} is, by \eqref{T9}.
Taking $p_1$ as variable, we infer from \eqref{T5}, \eqref{T7} and \eqref{T8}
that 
$$
\frac{1}{p_2}=\frac{1}{\alpha -1}(1-\frac{3}{p_1}),\quad
\frac{1}{q_1}=\frac{1}{3}(1+\frac{2}{n\alpha}) -\frac{2}{np_1},\quad
\frac{1}{q_2}=\frac{2}{n(\alpha -1)} (\frac{3}{p_1}-\frac{1}{\alpha}).
$$ 
The constraints \eqref{T4}, \eqref{T30} and \eqref{T40} are 
found to be respectively equivalent to 
\be
\label{constraints}
p_1\le 3(1-\frac{\alpha -1}{4})^{-1}  (\text{for }\ 
\alpha \le 4), \quad p_1\ge \sup \big\{ 6(n+\frac{2}{\alpha})^{-1},
6(1-\frac{2}{n})(1-\frac{4}{n\alpha})^{-1}\big\}, \quad
p_1 < 3\alpha .
\ee
The value $p_1=6$ fulfills all the requirements in 
\eqref{constraints}. Let now
$s>\frac{n}{2}-\frac{2}{\alpha }$ be given. Choose $\varepsilon >0$ such that
$4\varepsilon < s-(\frac{n}{2}-\frac{2}{\alpha})$, and pick
$s_1\in (\sigma _1, \sigma _1 + \varepsilon )$,  and 
$s_2\in (\sigma _2, \sigma _2 + \varepsilon )$. 
Then \eqref{S40} and \eqref{S41} hold for some numbers 
$b_1<\frac{1}{2}$, $b_2<\frac{1}{2}$, 
according to Lemma \ref{interpolation}. Set finally
$b=\sup \{b_1,b_2\}$. Then we have
$$
\Sigma \le c \big( \prod_{i=1}^{\alpha +1} ||u_i||_{X_{s,b}}\big) 
||v||_{X_{-s,b}}
$$ 
which gives \eqref{multilinear}.\\
\underline{(ii) $\alpha=2$}\\
Observe first that the approach followed in (i) does not work 
for $n>2$. Indeed,  the constraints \eqref{T3}-\eqref{S1000} 
impose $p_1=p_2=q_1=q_2=4$, and the equation $3\sigma _1=\sigma _2$ 
is then satisfied only for $n=2$. 
Assume $n\ge 3$. We now search a couple
$(p_1,q_1)$ satisfying
\be
0<\frac{1}{p_1} \le \frac{1}{q_1} \le \frac{1}{2} 
\le \frac{1}{p_1}+\frac{1}{q_1} \le 1,\quad
s_1> (n-2)(\frac{1}{2}-\frac{1}{q_1}), \quad
b_1> 1-\frac{1}{p_1}-\frac{1}{q_1},
\label{W1}
\ee
while $(p_2,q_2)$ still satisfies
\be
0 < \frac{1}{p_2}\le \frac{1}{4},\quad 
0\le \frac{1}{q_2}\le \frac{1}{2}-\frac{1}{p_2},\quad 
s_2>\frac{n}{2} -\frac{2}{p_2}-\frac{n}{q_2}\cdot 
\label{W2}
\ee
The H\"older exponents $(p_1,q_1)$ and $(p_2,q_2)$ have to satisfy 
the relations 
\begin{eqnarray}
\frac{3}{p_1}+\frac{1}{p_2} &=& 1, \label{W3}\\
\frac{3}{q_1}+\frac{1}{q_2} &=& 1. \label{W4}
\end{eqnarray}
We still minimize the functional $\sup \{ 3\sigma _1, \sigma _2\}$, where 
$$
\sigma _1=(n-2)(\frac{1}{2}-\frac{1}{q_1}), \quad
\sigma _2 = \frac{n}{2} -\frac{2}{p_2}-\frac{n}{q_2}=
\frac{n}{2}-\frac{2}{p_2} -n (1-\frac{3}{q_1})
$$ 
by solving in $q_1$ the equation
$3\sigma _1=\sigma _2$. Taking $p_2=4$ to
produce the least value of $\sigma _2$, we find as solution
$q_1=3(1+\frac{1}{4n-5})\in (3,4)$, which yields $p_1=4$ and 
$q_2=4(n-1)$ by \eqref{W3}-\eqref{W4}, and 
$$
3\sigma _1=\sigma _2= \frac{n}{2}-\frac{3}{4}-\frac{1}{4(n-1)}\cdot
$$
The constraints on 
$p_1,q_1,p_2,q_2$ in \eqref{W1}-\eqref{W2} are clearly fulfilled, 
for $n>2$. Pick now any $s>\frac{n}{2}-\frac{3}{4}-\frac{1}{4(n-1)}$
and $\varepsilon >0$ such that $4\varepsilon < s-(\frac{n}{2}-\frac{3}{4}
-\frac{1}{4(n-1)})$. 
We next pick 
$s_1\in (\sigma _1, \sigma _1 +\varepsilon )$,
$s_2\in (\sigma _2, \sigma _2 +\varepsilon )$, 
$b_1\in (1-\frac{1}{p_1}-\frac{1}{q_1},\frac{1}{2})$, and $b_2<\frac{1}{2}$
so that \eqref{S21} holds. 
Then  \eqref{multilinear} follows with $b=\sup\{ b_1,b_2\}$. \\
\underline{(iii) $\alpha =1$}\\
In this case, we have with $p_1=q_1=3$
$$
\Sigma \le c 
||u_1||_{X_{s,b_1}}
||u_2||_{X_{3s_1+\varepsilon ,b_1}}
||v||_{X_{-s,b_1}}
$$
provided that \eqref{W1} is satisfied, i.e.
$$
s_1>\sigma_1=\frac{n-2}{6},\quad b_1>\frac{1}{3}\cdot
$$ 
Therefore, if $s>\frac{n}{2}-1$, taking $\varepsilon >0$ such that
$4\varepsilon < s-(\frac{n}{2}-1)$, $s_1\in (\sigma _1, \sigma _1
+\varepsilon)$, and $b=b_1\in (\frac{1}{3},\frac{1}{2})$, we conclude that 
$$
\Sigma \le c||u_1||_{X_{s,b}} ||u_2||_{X_{s,b}} ||v||_{X_{-s,b}}
$$
and \eqref{multilinear} follows. \qed

\subsection{Proof of Proposition \ref{prop30}.}
We begin with the proof of \eqref{PP1} by following closely 
\cite{grunrock01}. Note, however, that the main concern here is to have 
the condition $s+2b<1/2$ fulfilled.  
Let $s,b$ be as in the statement of
Proposition \ref{prop30}, and let $v_1,v_2\in X_{s,b}$ be decomposed as 
$$
v_i(x,t)
=\int_\R \sum_{k\in \Z ^2} {\cal F} v_i (k, \tau ) 
e^{i(k\cdot x + \tau t)} d\tau \qquad i=1,2.
$$ 
(Here, we use the symbol $\cal F$ instead of $\hat \cdot$ to denote Fourier transform in space and time.)
Let 
$$
f_i(k, \tau ) = \la k\ra ^s \la \tau - |k|^2 \ra ^{b} 
{\cal F }\, {\overline{ v}_i}(k,\tau ),\quad \ i=1,2.
$$ 
Then 
\begin{equation}
 \label{PP3}
||\overline{v}_1\overline{v}_2||_{X_{s,b'}}
= || \la k\ra ^s \la \tau + |k|^2 \ra ^{b'}
\int_{\tau _1 + \tau _2 = \tau}  \sum_{k_1+k_2=k}\ 
\prod_{i=1}^2 \la k_i\ra ^{-s} \la \tau _i -|k_i|^2\ra ^{-b} f_i   
||_{L^2_{k,\tau}}
\end{equation}
where $\int_{\tau _1 + \tau _2 = \tau}  \sum_{k_1+k_2=k}$ stands for
$\int_{\R} d\tau _1 \sum_{k_1\in \Z ^2}$ with the relations  
$\tau _1 + \tau _2 = \tau$ and  $k_1+k_2=k$ satisfied. 
Let $A_0$ (resp.  $A_i$, $i=1,2$) denote the region where the largest number 
among $\la \tau +|k|^2\ra$, $\la\tau _1 -|k_1|^2 \ra $ and 
$\la \tau _2 -|k_2|^2 \ra $, is $\la \tau +|k|^2 \ra $
(resp. $\la \tau _i -|k_i|^2\ra$, $i=1,2$). We infer from the relation
$$
\tau + |k|^2 -\sum_{i=1}^2 (\tau _i - |k_i|^2) 
= |k|^2 +\sum_{i=1}^2  |k_i|^2
$$ 
that 
\begin{equation}
\label{PP4}
\la k\ra ^2 +\sum_{i=1}^2 \la k_i \ra ^2 \le 
C\left( \la \tau + |k|^2 \ra  + \sum_{i=1}^2 \la \tau _i -|k_i|^2 \ra 
\right) 
\end{equation}
Let us begin with the region $A_0$. 
\eqref{PP4}  gives, with $0<\varepsilon < 
\inf \{ \frac{1}{2}(\frac{1}{2}-|s|), 2(b-|s|) \} $
and $-b':=\frac{1}{2}(\frac{1}{2}-s) + \varepsilon <\frac{1}{2}$
$$
\la k\ra ^{\frac{1}{2} +s} \prod_{i=1}^2 \la k_i\ra ^{-s+\varepsilon}
\le C \la \tau + |k|^2 \ra ^{-b'}.
$$
The contribution in \eqref{PP3} due to $A_0$ is therefore bounded by 
\begin{eqnarray*}
&&C||\la k\ra ^{-\frac{1}{2}} 
\int_{\tau _1 + \tau _2=\tau } \sum_{k_1+k_2=k} 
\la k_i\ra ^{-\varepsilon} \la \tau _i -|k_i|^2\ra ^{-b}
|f_i| ||_{L^2_{k,\tau}} \\
&&\qquad =C ||\la k\ra ^{-\frac{1}{2}}
\int_{\tau _1 + \tau _2 =\tau}\sum_{k_1+k_2=k}
\la k_i\ra ^{s-\varepsilon } |{\cal F}\, \overline{v}_i| ||_{L^2_{k, \tau}} \\
&&\qquad  =C||\prod_{i=1}^2 J^{s-\varepsilon} 
{\cal F}^{-1} |{\cal F}\, {\overline v_i} | 
||_{L^2_tH^{-\frac{1}{2}}_x} \\
&&\qquad  \le C||\prod_{i=1}^2 J^{s-\varepsilon} 
{\cal F}^{-1} |{\cal F}\, {\overline v_i} | 
||_{L^2_tL^q_x}, \qquad q>\frac{4}{3}\\
&&\qquad \le C\prod_{i=1}^2 ||J^{s-\varepsilon}
{\cal F}^{-1} |{\cal F}\, {\overline v_i} | 
||_{L^4_tL^{2q}_x}, \qquad q>\frac{4}{3}\\
&&\qquad  \le C\prod_{i=1}^2 ||J^{s-\varepsilon} 
{\cal F}^{-1} |{\cal F}\, {\overline v_i} | 
||_{X^-_{\varepsilon, b}} \\
&&\qquad  \le C\prod_{i=1}^2 ||v_i||_{X_{s,b}}
\end{eqnarray*}
where we used the fact that $L^q(\T ^2 ) \subset H^{-\frac{1}{2}}(\T ^2)$ 
for $q>4/3$ (by dualizing the Sobolev embedding $H^\frac{1}{2}(\T ^2) \subset L^p(\T ^2)$ for $p<4$), H\"older inequality, and \eqref{S21}-\eqref{S22}.
We also used the notation 
$$||u||_{X^-_{s,b}}=(\int_\R \sum_{k\in \Z ^2} \la k\ra ^{2s}
\la \tau -|k|^2\ra ^{2b} |{\cal F}u(k,\tau )|^2 d\tau )^{\frac{1}{2}} =
||\overline{u}||_{X_{s,b}}$$ 
borrowed from \cite{grunrock00}.
It remains to estimate the contributions in \eqref{PP3} due to the regions 
$A_1$ and $A_2$. By symmetry, we can consider only the region $A_1$. 
In $A_1$, since $-s+\frac{\varepsilon}{2} <b$, we have that
$$
\la k_2\ra ^{-s+\varepsilon} \la k_1\ra ^{-s} 
\le C \la \tau _1 - |k_1|^2 \ra ^{-s+\frac{\varepsilon}{2}}
\le C \la \tau _1 -|k_1|^2\ra ^b 
$$ 
and therefore the contribution in \eqref{PP3} is bounded by 
$$
||\la k\ra ^s \la \tau + |k|^2 \ra ^{b'} 
\int_{\tau _1 + \tau _2 =\tau }\sum_{k_1+k_2 =k} 
|f_1|\la k_2\ra ^{-\varepsilon}
\la \tau _2 -|k_2|^2\ra ^{-b} |f_2| ||_{L^2_{k,\tau}}
=C||{\cal F}^{-1} |f_1| J^{s-\varepsilon} {\cal F}^{-1}
|{\cal F}\, \overline{v}_2 |  ||_{X_{s,b'}}.
$$
By \eqref{S20}-\eqref{S21} with $-s>1/3$ and $-b'$ chosen sufficiently close to 
$\frac{1}{2}$, we have that
$$
X_{-s,-b'}\subset L^6(\R ; L^6(\T ^2)),\quad \text{ hence }\quad 
L^{\frac{6}{5}}(\R ;  L^{\frac{6}{5}}(\T ^2)) \subset X_{s,b'}.
$$
It follows that 
\begin{eqnarray*}
||{\cal F}^{-1} |f_1| J^{s-\varepsilon} {\cal F}^{-1}
|{\cal F}\, \overline{v}_2 |  ||_{X_{s,b'}}
 &\le& 
C||{\cal F}^{-1} |f_1| J^{s-\varepsilon}{\cal F}^{-1}
|{\cal F}\, \overline{v} _2| ||_{L^{\frac{6}{5}}_tL^{\frac{6}{5}}_x}\\
&\le& C||{\cal F}^{-1} |f_1|||_{L^2_tL^2_x}
|| J^{s-\varepsilon}{\cal F}^{-1}
|{\cal F}\, \overline{v}_2| ||_{L^3_t L^3_x}\\
&\le& C||\overline{v}_1||_{X^-_{s,b}}
|| J^{s-\varepsilon }{\cal F}^{-1}|{\cal F}\, \overline{v}_2| 
||_{X^-_{\varepsilon ,b}}\\
&\le& C||v_1||_{X_{s,b}}
|| v_2 ||_{X_{s,b}}\\
\end{eqnarray*}
where we used H\"older inequality and \eqref{S21}-\eqref{S22} with $p=q=3$.
This completes the proof of \eqref{PP1}. 

To derive \eqref{PP2} from \eqref{PP1}, we consider two functions $u_1,u_2$
in $X_{0,b}(\Omega ) \subset X_{s,b}(\Omega )$, and consider their odd extensions
$v_1,v_2$ to $(-\pi ,\pi)^2$; i.e., $v_i(\epsilon _1 x_1,\epsilon _2 x_2)
=\epsilon_1\epsilon_2 u_i(x_1,x_2)$ for $x=(x_1,x_2)\in \Omega$ and 
$\epsilon_i =\pm 1$. Note that $v_1,v_2\in X_{0,b}$ and that 
$\overline{u}_1\overline{u}_2=(\overline{v}_1\overline{v}_2)_{\vert _{\Omega }}$.
For any function 
$w=\sum_{k\in \N ^2} \int_\R {\cal F} w(k, \tau)
e^{i\tau t}\cos (k_1x_1 ) \cos (k_2x_2) d\tau$, we set 
$$
||w||^2_{X_{s,b}(\Omega )_N} =
\sum_{k\in \N ^2} \int_\R \la \tau + |k|^2 \ra ^{2b}
\la k\ra ^{2s} |{\cal F} w(k,\tau )|^2 d\tau.
$$
The Bourgain space $X_{s,b}(\Omega )_N$ (with Neumann boundary conditions) is 
defined as the space of the $w$'s for which the norm 
$||w||_{X_{s,b}(\Omega )_N}$ is finite. Since the function $\overline{v}_1\overline{v}_2$ is even with respect to both $x_1$ and $x_2$, we have that 
$$
||\overline{u}_1\overline{u}_2||_{X_{s,b'}(\Omega )_N} 
\sim C ||\overline{v}_1\overline{v}_2 ||_{X_{s,b'}} 
\le C||v_1||_{X_{s,b}} ||v_2||_{X_{s,b}}
\le C ||u_1||_{X_{s,b}(\Omega )} ||u_2||_{X_{s,b}(\Omega )}\cdot
$$
We claim that $X_{s,b}(\Omega ) =X_{s,b}(\Omega )_N$ for $|s|<1/2$ and 
$|b|\le 1$. 
Note first that this is true for $|s|<\frac{1}{2}$ and $b=0$, since 
$$
X_{s,0}(\Omega )=L^2(\R ; H^s(\Omega ))=X_{s,0}(\Omega )_N.
$$
The claim is also true for $|s|<1/2$ and $b=1$, since 
$$
u\in X_{s,1}(\Omega ) \iff u\in X_{s,0}(\Omega ) \text{ and } 
iu_t+\Delta u \in X_{s,0}(\Omega ) 
$$
and since a similar criterion may be written for $X_{s,1}(\Omega )_N$.
The claim is also true for $|s|<1/2$ and $0\le b\le 1$ by interpolation,
and for $|s|<1/2$ and $|b|\le 1$ by duality. \eqref{PP2} follows for
$u_1,u_2\in X_{0,b}(\Omega )$, and also for 
$u_1,u_2\in X_{s,b}(\Omega )$ by density.
This completes the proof of Proposition \ref{prop30}.

\end{document}